\documentclass[preprint,12pt, a4paper]{elsarticle}

\usepackage{amssymb}
\usepackage{amsmath}
\usepackage{mathtools}
\usepackage[dvipsnames]{xcolor}
\usepackage{placeins}
\usepackage{multirow}
\usepackage{hyperref}
\usepackage{subcaption}
\usepackage{verbatim}

\usepackage[a4paper,margin=2cm]{geometry}

\journal{}

\begin{document}

\begin{frontmatter}



\title{Development of a GPU-accelerated numerical simulation model for liver infection in realistic domains}



\author[label1]{Eva Cirujeda}
\ead{873283@unizar.es}

\author[label1]{Adrián Navas-Montilla}
\ead{anavas@unizar.es}

\affiliation[label1]{organization={Fluid Dynamics Technologies Group, Aragon Institute of Engineering Research (I3A), University of Zaragoza},
            addressline={}, 
            city={Zaragoza},
            postcode={}, 
            country={Spain}}

\author[label2]{Cordula Reisch}
\ead{cordular@ruc.dk}

\affiliation[label2]{organization={Department of Science and Environment, Roskilde University},
            addressline={Universitetsvej 1}, 
            city={Roskilde},
            postcode={4000}, 
            country={Denmark}}

\date{November 2026}

\begin{abstract}
Viral hepatitis reduces the functionality of the liver and constitutes a major global health concern because it oftentimes evolves to chronic inflammation leading to long term liver damage. Mathematical models on a mesoscopic scale abstract from only partly known cell-scale processes while providing an in-silico tool for understanding tissue-scale development of the disease. In this work, we investigate the interactions between virus and immune cells in a realistic three-dimensional liver geometries by using a reaction-diffusion model that reproduces infection courses with cure and with chronic developments. A GPU-based finite-difference implementation, released as an open-source Python library, achieves speedups greater than $\times 30$ for realistic mesh sizes, making three-dimensional simulations on realistic geometries computationally tractable. This provides a basis for future in-silico studies of treatment strategies, such as testing antiviral or immunomodulatory therapies, and for patient-specific simulations based on individual liver geometries, paving the way towards personalized healthcare.

\end{abstract}

\begin{keyword}
partial differential equations \sep numerical modelling \sep HPC computing \sep finite differences \sep mathematical biology
\end{keyword}

\end{frontmatter}

\section{Introduction}
\label{sec:introduction}
\begin{figure}[!h]
    \centering
    \includegraphics[width=1.0\textwidth]{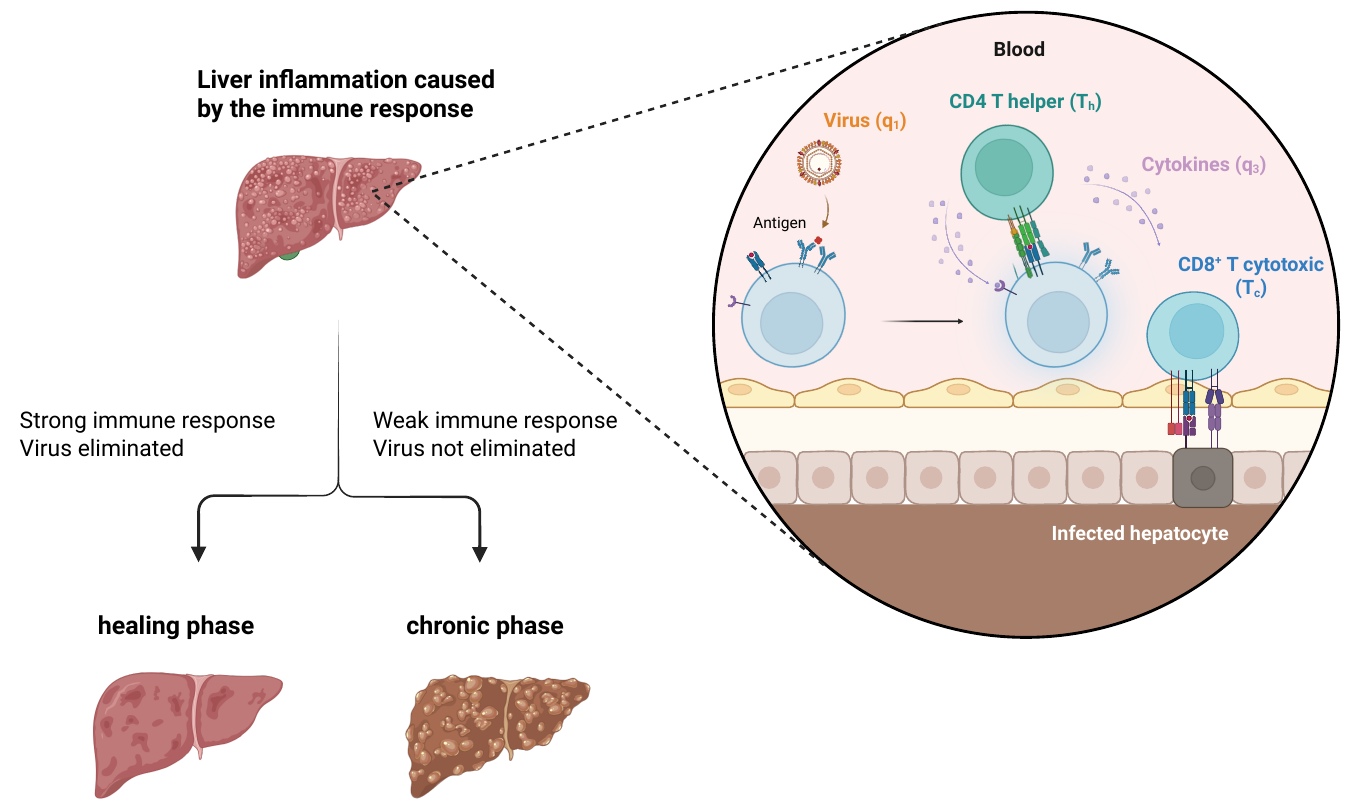}
    \caption{Viral infection pathway, created with BioRender (https://BioRender.com).}
    \label{fig:Ruta de infección del virus}
\end{figure}

Hepatitis B is an infectious liver disease caused by the hepatitis B virus (HBV), which specifically infects hepatocytes, the functional cells of liver tissue. HBV infection is a major global public health problem because it can be either acute or chronic; chronic infection increases the risk of long-term complications such as liver cirrhosis and hepatocellular carcinoma \cite{Usuda2024}. Liver damage is not caused directly by the virus, but by the adaptive immune response triggered after infected hepatocytes are detected. In this process, helper T lymphocytes ($T_h$) coordinate the response through the secretion of cytokines, which activate cytotoxic T lymphocytes ($T_c$), enhance their destructive activity, and promote their migration towards infected regions. Although this immune activity aims to reduce the viral load, it also causes cellular damage that may give rise to fibrotic tissue. In turn, this tissue hinders immune-cell transport and may promote the chronification of the infection rather than its resolution. The Figure \ref{fig:Ruta de infección del virus} schematically summarises the route of viral infection, the associated immune response and the two possible scenarios considered in this work (virus elimination and chronification of the infection).

Mathematical modelling offers a perspective complementary to clinical and experimental studies, as it makes it possible to abstract the still largely unknown infection mechanisms at the cellular scale while reproducing the population-scale interaction between the virus and the immune response. One of the first models along these lines was proposed by Kerl et al. \cite{Kerl2012}, who used a predator--prey reaction--diffusion system to show that certain geometric features of the domain, such as internal obstacles, can drive the infection towards a chronic steady state. Reisch and Langemann \cite{Reisch2019} subsequently extended this model by incorporating a chemotactic advection term to describe the directed migration of immune cells along chemical gradients. They concluded that, although advection alone does not explain chronification, it improves the physical realism of cell-transport representation. On this basis, Reisch et al. \cite{Reisch2024} introduced a more general family of reaction--diffusion models for liver inflammation that explicitly included the virus, helper and cytotoxic T lymphocytes, and cytokines. The authors also noted the difficulty of validating these models against real data because of the scarcity of spatiotemporal clinical measurements. More recently, Foko and Mohammadi \cite{Foko2026} proposed a diffusive model with time delays, viral capsids, and mitotic transmission, solved in three-dimensional domains. This work reflects the trend towards increasingly detailed spatial descriptions of hepatitis B dynamics. On the other hand,  space-independent models are studied with a focus on more detailed descriptions of the immune mechanisms, \cite{BoivinChampeaux2026}. 

As the dimensionality and spatial resolution of these models increase, so does their computational cost because of the larger number of mesh nodes, the additional operations required at each time step, and the volume of data that must be stored during the simulation. Graphics processing units (GPUs) are an effective and widely used alternative for reducing this cost, as they allow the workload to be distributed across thousands of parallel processing cores. GPU acceleration has proven effective in several biomedical contexts, including nonlinear finite-element simulations for surgical applications \cite{Taylor2008}, three-dimensional models of cell growth and division \cite{Madhikar2018}, and agent-based biological simulations \cite{Hesam2021}, commonly reducing computation times by an order of magnitude.

Despite these advances, reaction--diffusion models developed for hepatitis B have so far been analysed mainly in non-spatial, one- and two-dimensional formulations and, to the best of our knowledge, have not yet been extended to anatomically realistic three-dimensional liver geometries or accelerated using GPU computing. In this work, we address this gap by: (i) extending the model of Reisch et al. \cite{Reisch2024} to three spatial dimensions; (ii) developing CPU and GPU implementations in Python using NumPy and CuPy of a second-order finite-difference scheme, achieving speedups greater than $\times 30$ for realistic mesh sizes; (iii) applying the resulting model to a three-dimensional liver geometry segmented from real computed-tomography data, including the vascular structure as the entry point for the immune response; and (iv) releasing the implementation as an open-source Python library installable through \texttt{pip} \textit{\href{https://github.com/navasmontilla/HepGPU}{HepGPU}} to facilitate reproducibility. Taken together, these contributions overcome the limitations of idealised domains and make it possible to investigate how spatial transport and tissue heterogeneity, particularly the presence of internal barriers and fibrotic regions, influence viral clearance or infection chronification.

\section{Mathematical model}
\label{sec:model}
\subsection{Definition of the model}
This section presents the mathematical model used to describe hepatitis B infection dynamics through a system of reaction--diffusion equations with chemotaxis. The model originates from \cite{Reisch2024} and describes the interaction among the virus density $q_1(\mathbf{x},t)$; helper and cytotoxic T lymphocyte densities $T_h(\mathbf{x},t)$ and $T_c(\mathbf{x},t)$, and the cytokine density $q_3(\mathbf{x},t)$.

\begingroup
\footnotesize
\begin{align}
  q_{1,t}
  &=
  \underbrace{a_1 q_1 (C_1 - q_1)
  \frac{q_1 - \varepsilon}{q_1 + \kappa}
  - a_5 q_1 T_c}_{\text{reaction term}}
  + \underbrace{\nabla \cdot
  \left[
  d_1^{\mathrm{ctc}}(\mathbf{x}) \nabla q_1
  \right]}_{\text{diffusion term}},
  \label{eq:q1}
  \\[1ex]
  T_{h,t}
  &=
  \underbrace{a_{2,h}\,\chi_{\Theta}(\mathbf{x})
  \left(C_{T_h} - T_h\right)
  \int_{\Omega} q_1 \, \mathrm{d}\mathbf{x}
  - a_{6,h} T_h (C_1 - q_1)}_{\text{reaction term}}
  + \underbrace{\nabla \cdot
  \left[
  d_{T_h}^{\mathrm{diff}}(\mathbf{x}) \nabla T_h
  \right]}_{\text{diffusion term}},
  \label{eq:Th}
  \\[1ex]
  T_{c,t}
  &=
  \underbrace{a_{2,c}\,\chi_{\Theta}(\mathbf{x})
  \left(C_{T_c} - T_c\right)
  \int_{\Omega} q_1 \, \mathrm{d}\mathbf{x}
  - a_{6,c} T_c (C_1 - q_1)}_{\text{reaction term}}
  - \underbrace{\nabla \cdot
  \left[
  d_{T_c}^{\mathrm{chem}}(\mathbf{x}) T_c \nabla q_3
  \right]}_{\text{chemotaxis term}},
  \label{eq:Tc}
  \\[1ex]
  q_{3,t}
  &=
  \underbrace{a_3 T_h q_1
  - a_{nd} q_3}_{\text{reaction term}}
  + \underbrace{\nabla \cdot
  \left[
  d_3^{\mathrm{diff}}(\mathbf{x}) \nabla q_3
  \right]}_{\text{diffusion term}},
  \label{eq:q3}
\end{align} 
\endgroup
where $\mathbf{x}=(x,y,z)$ denotes spatial position, $t\geq 0$ denotes time, and nabla, $\nabla=(\frac{\partial}{\partial x}, \frac{\partial}{\partial y}, \frac{\partial}{\partial z})$, is the vector differential operator. The model is defined on the spatial domain $\Omega \subset \mathbb{R}^3$, which has subdomains $\Theta \subset \Omega$.

\subsection{Modelled mechanisms}

Two spatial propagation mechanisms are considered in the model: First, the diffusion term describes the spatial propagation of the densities throughout the domain. It models flow from regions of higher densities towards regions of lower densities, thereby tending to homogenise or smooth spatial gradients. On the other hand, the chemotaxis term appears only in the equation for cytotoxic T lymphocytes $T_c$ and describes transport driven by the cytokine $q_3$ gradient; cytotoxic T lymphocytes move towards regions with higher cytokine concentrations. This term will be treated numerically as an advection term.

The reaction term models the local dynamics of each of the four model variables, describing their growth and decay. For the virus $q_1$, growth depends on two mechanisms. First, logistic growth term $a_1 q_1 (C_1 - q_1)$ limits viral proliferation as it approaches the environmental carrying capacity $C_1$. Second, the Allee effect $\frac{q_1 - \varepsilon}{q_1 + \kappa}$ reduces the growth rate at low viral densities. The denominator scales the Allee effect closer to a pure logistic growth for densities close to the carrying capacity.
Viral decay is related to the elimination rate $a_5$ exerted by cytotoxic T lymphocytes $T_c$.

To describe the two T-lymphocyte variables, the generic variable $q_2 = (T_h, T_c)$ is introduced. Its growth is characterised by a logistic term,
\begin{align}
  a_{2,i}\,\chi_{\Theta}(\mathbf{x})\left(C_{T_i} - q_2\right)
  \int_{\Omega} q_1 \, \mathrm{d}\mathbf{x},
  \quad i \in \{h,c\},
\end{align}
where the subdomain $\Theta$ represents the T-lymphocyte entry region. 
The function $\chi_{\Theta}(\mathbf{x})$ represents a localised inlet associated with the vascular region and models the flow of T lymphocytes from the bloodstream into the liver domain. $\xi(\mathbf{x})$ is defined from a mask specifying the size and location of the domain inlet; it is nonzero only in this region and zero elsewhere in the domain.
\begin{align}
\xi(\mathbf{x})=
\begin{cases}
1, & \mathbf{x}\in \Theta,\\
0, & \mathbf{x}\in \Omega\setminus \Theta,
\end{cases}
\label{eq:chiTheta_x}
\end{align} 
The function $\xi(\mathbf{x})$ is normalised:
\begin{align}
\chi_{\Theta}(\mathbf{x}) = \frac{\xi(\mathbf{x})}{\int_{\Omega} \xi(\mathbf{x})\, d\mathbf{x}},
\label{eq:chiTheta_x_normalizada}
\end{align} 
so that its total contribution over the domain is unity,
\begin{align}
\int_{\Omega}\chi_{\Theta}(\mathbf{x})\, d\mathbf{x}=1.
\label{eq:chiTheta_dominio_unitario}
\end{align} 

Activation of the immune response through these variables depends on the total viral load in the domain ($\int_{\Omega} q_1 \, \mathrm{d}\mathbf{x}$), while growth slows as the maximum capacity $C_{T_i}$ is approached. This mechanism describes the interaction with the immune system of the whole body because the production of immune cells takes place outside of the liver tissue and depends on the viral load within the whole organ. 
The natural degradation of the immune cells, modelled by $a_{6,i},q_2,(C_1 - q_1)$, regulates their disappearance according to the absence of virus.

Finally, for the cytokines $q_3$, the reaction term is given by $a_3 T_h q_1 - a_{nd} q_3$. The first term represents cytokine production induced by the interaction between helper T lymphocytes $T_h$ and the virus, whereas the second represents their natural degradation in the medium.

\subsection{Model parameters}

The model parameters are detailed in the corresponding tables. Table \ref{tab:parametros_locales_modelo} describes the local parameters, which have no spatial dependence, whereas the parameters associated with the spatially dependent diffusion and advection terms are given in Table \ref{tab:parametros_distribuidos_modelo}.
The parameter values are artificial, not representing biological mechanisms in reality. 

\begin{table}[htbp]
\centering
\caption{Model parameters of the reaction terms}
\label{tab:parametros_locales_modelo}
\resizebox{\textwidth}{!}{
\begin{tabular}{l p{11cm} c}
\hline
\textbf{Parameter} &
\textbf{Description} &
\textbf{Value} \\
\hline

$a_{1}$ &
Growth or natural replication rate of the virus $q_1$ &
1.0 \\

$a_{3}$ &
Production rate of cytokines $q_3$ &
0.8 \\

$a_{nd}$ &
Disintegration rate of cytokines $q_3$ &
0.6 \\

$a_{5}$ &
Effectiveness of cytotoxic T cells $T_c$ for viral elimination &
0.05 \\

$a_{2,h}$, $a_{2,c}$ &
Regulation of the entry of T cells ($T_h$, $T_c$) through $\Theta$ &
2.0 \\

$a_{6,h}$, $a_{6,c}$ &
Decay rate of T cells ($T_h$, $T_c$) &
0.2 \\

$C_{1}$ &
Maximum capacity of the virus $q_1$ &
1.0 \\

$C_{Th}$ &
Maximum capacity of helper T cells $T_h$ &
8.0 \\

$C_{Tc}$ &
Maximum capacity of cytotoxic T cells $T_c$ &
15.0 \\

$\epsilon$ &
minimal Allee capacity &
0.05 \\

$\kappa$ &
Allee effect parameter &
0.01 \\

\hline
\end{tabular}}
\end{table}
\FloatBarrier

\begin{table}[htbp]
\centering
\caption{Diffusion parameters of the model (with spatial dependence)}
\label{tab:parametros_distribuidos_modelo}
\resizebox{\textwidth}{!}{
\begin{tabular}{l p{11cm} c}
\hline
\textbf{Parameter} &
\textbf{Description} &
\textbf{Value} \\
\hline

$d_1^{\mathrm{ctc}}$ &
Diffusion coefficient of the virus $q_1$ &
0.6 \\

$d_{T_h}^{\mathrm{diff}}$ &
Diffusion coefficient of helper T lymphocytes $T_h$ &
0.9 \\

$d_{T_c}^{\mathrm{chem}}$ &
Advection coefficient of cytotoxic T lymphocytes $T_c$ &
0.5 \\

$d_3^{\mathrm{diff}}$ &
Diffusion coefficient of cytokines $q_3$ &
0.5 \\

\hline
\end{tabular}}
\end{table}
\FloatBarrier

\section{Numerical method and computational implementation} 
\label{sec:metodo_discretizacion}

This section presents the discretisation method in a single spatial dimension, $x$, since its extension to higher dimensions is performed analogously by applying the same procedure in each direction. The spatial domain is taken to be one-dimensional, $x \in [0, L]$, while the temporal evolution of the system is studied over the interval $t \in [0, t_f]$. Spatial nodes $x_i$ and time levels $t^n$ are introduced for the discretisation, so that the solution is approximated by discrete values $u_i^n \approx u(x_i, t^n)$.

The spatial domain is discretised using a uniform mesh of $N$ nodes at positions $x_i$, which defines the spatial step as
\begin{align}
\Delta x = \frac{L}{(N-1)}.
\label{paso_espacial_dominio_x}
\end{align}
For the temporal discretisation, the time interval is divided into discrete levels $t^n$, where the index $n$ denotes the corresponding time and the time step between two consecutive times is $\Delta t = t^n-t^{n-1}$.

The system of equations defining the mathematical model is given in Equations \eqref{eq:q1}, \eqref{eq:Th}, \eqref{eq:Tc}, and \eqref{eq:q3}. From a numerical perspective, two types of equations are identified. The variables \(q_1\), \(T_h\), and \(q_3\) follow reaction--diffusion dynamics, with the spatial transport modelled by a second-order derivative:
\begin{align}
\frac{\partial u (x,t)}{\partial t}
&=
\frac{\partial }{\partial x}\left(\alpha(x)\frac{\partial u (x,t)}{\partial x}\right)
+ R(u),
\label{eq:generica_dif}
\end{align} where $u$ may denote the viral load (\(q_1\)), helper T lymphocytes (\(T_h\)), or cytokines (\(q_3\)), and \(\alpha\) is the diffusion coefficient governing the spatial propagation of $u$ throughout the domain. The diffusion term is given by \(\frac{\partial }{\partial x}\left(\alpha(x)\frac{\partial u (x,t)}{\partial x}\right)\), while \(R(u)\) represents the reaction terms describing local biological interactions among the system variables.

The variable \(T_c\), in contrast, follows advection--reaction dynamics characterized by a first-order spatial derivative:
\begin{align}
\frac{\partial u(x,t)}{\partial t}
=
- \frac{\partial f(u) }{\partial x}
+ R(u),
\label{eq:generica_adv}
\end{align}
where $u$ represents cytotoxic T cells ($T_c$), $f(u)=au(x,t)$ is the transport flux, and \(a\) is the transport velocity. The advection term is expressed as \(\frac{\partial f(u)}{\partial x}\) and the term \(R(u)\) represents the reaction terms describing local biological interactions among the system variables.

The reaction terms are evaluates explicitly at the mesh nodes at the current time, $R(u_i^n)$.

\subsection{Discretisation of the diffusion terms}
For the diffusion terms, the second spatial derivative is approximated using central finite differences. This approximation is obtained from the Taylor-series expansion of the function about a mesh point $u_i^n$.
\begin{align}
  \left.\dfrac{\partial}{\partial x}\left(\alpha\dfrac{\partial u}{\partial x}\right)\right|_{x_i}
  \;\approx\;
  \dfrac{\widetilde{\alpha}_{i+1/2}(u_{i+1}^n - u_i^n) - \widetilde{\alpha}_{i-1/2}(u_{i}^n - u_{i-1}^n)}{\Delta x^2}
\end{align}
where $\widetilde{\alpha}_{i+1/2}$ and $\widetilde{\alpha}_{i-1/2}$ are the coefficients evaluated between two mesh nodes and are defined as
\begin{align}
  \widetilde{\alpha}_{i+1/2} =
  \begin{cases}
  \dfrac{\alpha_{i+1} + \alpha_i}{2}, & \text{if } \alpha_{i+1} > 0 \text{ and } \alpha_i > 0, \\
  0, & \text{otherwise},
  \end{cases}
\end{align}
and analogously for $\widetilde{\alpha}_{i-1/2}$. It is noted that the central differences approximation  has an error of order \(O(\Delta x^2)\).

\subsection{Discretisation of the advection term}
The advection term is discretised using a first order upwind scheme. This term is approximated by using the definition of numerical fluxes between nodes, at $(i+\tfrac{1}{2})$ and $(i-\tfrac{1}{2})$, as follows:
\begin{align}
  \frac{\partial f(u)}{\partial x}
  \;\approx\;
  \frac{f_{i+1/2}-f_{i-1/2}}{\Delta x},
\end{align}
where $f_{i+1/2}$ and $f_{i-1/2}$ are the numerical fluxes entering and leaving the node. The fluxes are computed using an upwind strategy in which the variable value is taken from the upstream node according to the direction of the velocity $\widetilde{a}$. In particular, the flux at the intermediate point $(i+\tfrac{1}{2})$ is defined as shown in Equation (\ref{eq:flujo_numerico_entrante}), with an analogous expression for $f_{i-1/2}$, as described in \cite{Hirsch2007},
\begin{align} 
  f_{i+1/2}
  =
  \begin{cases}
  \widetilde{a}_{i+1/2}\,u_i^n, & \text{if  }  \widetilde{a}_{i+1/2}>  0  \\[8pt]
  \widetilde{a}_{i+1/2}\,u_{i+1}^n, & \text{if  } \widetilde{a}_{i+1/2} < 0 
  \end{cases}
  \label{eq:flujo_numerico_entrante}
\end{align}
where $\widetilde{a}_{i+1/2}$ is the approximated advection velocity at the intermediate point $(i+\tfrac{1}{2})$ between two nodes, obtained from the gradient of the variable driving transport. In the case considered here, this velocity is defined as a function of the cytokine gradient $\nabla q_3$:
\begin{align}
  \widetilde{a}_{i+1/2} = d_{Tc,i+1/2}\,\frac{q_{3,i+1}^n - q_{3,i}^n}{\Delta x},
\end{align}
where $d_{Tc,i+1/2}$ denotes the transport coefficient at the intermediate point.

\subsection{Discretisation of the temporal derivative}
The temporal discretisation is performed using the first-order explicit Euler scheme.
To ensure numerical stability, the dynamic time step $\Delta t$ is calculated at each
iteration according to the advective and diffusive constraints. Specifically, the time
step is limited by the Courant--Friedrichs--Lewy (CFL) condition for the advective terms
and by the more restrictive condition associated with the second-order derivative for
the diffusive terms:
\begin{equation}\label{eq:sigma}
\Delta t_k = \sigma \frac{\Delta x^2}{\max_{i}(\tilde{\alpha}_{k,i+1/2})}, \quad 
\Delta t_{Tc} = \sigma \frac{\Delta x}{\max_{i}(\tilde{a}_{i+1/2})},
\end{equation}
where $\sigma \leq 0.5$ and $\sigma \leq 1$ act as safety factors, respectively. The final
global time step used in the simulation is determined adaptively as
\begin{equation}
\Delta t = \min(\Delta t_1, \Delta t_{Th}, \Delta t_{Tc}, \Delta t_3)
\end{equation}
thereby ensuring that all stability criteria are satisfied simultaneously.

\subsection{Boundary conditions}
The preceding expressions update the interior nodes of the domain, but boundary nodes are treated differently. Specifically, a homogeneous Neumann condition (zero normal flux) is imposed on all four model variables along the entire outer boundary, irrespective of the spatial dimensionality. The same boundary condition is applied at the edge of every internal obstacle (barriers or fibrotic regions) by setting the model variables and the diffusion and advection coefficients to zero within those obstacles, consistently with the discretisation of the transport operators described in Section \ref{sec:metodo_discretizacion}.

\section{Numerical Validation}
\label{sec:numerical_validation}

\subsection{Validation of the space-time convergence order}

The spatial convergence of the scheme is evaluated in the one-dimensional setting through a numerical mesh-refinement study. 
Solutions were computed for five consecutive mesh sizes, $N=\{11,21,41,81,161\}$ nodes, up to $t_f=20$ in a domain of length $L=1$, with the remaining parameters set according to Table \ref{tab:parametros_locales_modelo}. For each pair of consecutive meshes, the coarse-mesh solution was compared with the fine-mesh solution evaluated at coincident nodes, and the differences were quantified using the $L_\infty^e$ and $L_2^e$ norms,
\begin{equation}
    L_{\infty}^e(u)
    =
    \max_{x \in [0,L]}
    \left| u_{\mathrm{coarse}}(x) - u_{\mathrm{fine}}(x) \right|,
    \label{eq:norma_linf_validacion_1D}
\end{equation}
\begin{equation}
    L_2^e (u)
    =
    \left( \sum_{i} \left| u_{\mathrm{coarse}}(x_i) - u_{\mathrm{fine}}(x_i) \right|^2 \, \Delta x \right)^{1/2}.
    \label{eq:norma_l2_validacion_1D}
\end{equation}
The results obtained for the four variables are reported in Tables \ref{tab:conv_Linf} and \ref{tab:conv_L2}.
\begin{table}[h!]
\centering
\caption{Spatial convergence errors in the $L_\infty^e$ norm}
\label{tab:conv_Linf}
    \resizebox{\textwidth}{!}{
    \begin{tabular}{c c c c c c c c}
    \hline
    $N_{\text{coarse}}$ & $N_{\text{fine}}$ & $\Delta x_{\text{c}}$ & $\Delta x_{\text{f}}$ &
    $L_\infty^e(q_1)$ & $L_\infty^e(T_h)$ & $L_\infty^e(T_c)$ & $L_\infty^e(q_3)$ \\
    \hline
    11  & 21  & 0.1 & 0.05  & $6.14\times10^{-2}$ & $1.52\times10^{-1}$ & $2.23$ & $7.55\times10^{-1}$ \\
    21  & 41  & 0.05 & 0.025  & $4.13\times10^{-2}$ & $9.07\times10^{-2}$ & $1.71$ & $4.93\times10^{-1}$ \\
    41  & 81  & 0.025 & 0.0125  & $2.45\times10^{-2}$ & $5.00\times10^{-2}$ & $1.31$ & $2.86\times10^{-1}$ \\
    81  & 161 & 0.0125 & 0.00625  & $1.34\times10^{-2}$ & $2.63\times10^{-2}$ & $1.01$ & $1.54\times10^{-1}$ \\
    \hline
    \end{tabular}}
\end{table}
\FloatBarrier
\begin{table}[h!]
\centering
\caption{Spatial convergence errors in the $L_2^e$ norm}
\label{tab:conv_L2}
    \resizebox{\textwidth}{!}{
    \begin{tabular}{c c c c c c c c}
    \hline
    $N_{\text{coarse}}$ & $N_{\text{fine}}$ & $\Delta x_{\text{c}}$ & $\Delta x_{\text{f}}$ &
    $L_2^e(q_1)$ & $L_2^e(T_h)$ & $L_2^e(T_c)$ & $L_2^e(q_3)$ \\
    \hline
    11  & 21  & 0.1 & 0.05  & $6.04\times10^{-2}$ & $1.25\times10^{-1}$ & $1.28$ & $7.54\times10^{-1}$ \\
    21  & 41  & 0.05 & 0.025  & $3.97\times10^{-2}$ & $7.52\times10^{-2}$ & $8.86\times10^{-1}$ & $4.91\times10^{-1}$ \\
    41  & 81  & 0.025 & 0.0125  & $2.33\times10^{-2}$ & $4.17\times10^{-2}$ & $5.60\times10^{-1}$ & $2.84\times10^{-1}$ \\
    81  & 161 & 0.0125 & 0.00625  & $1.27\times10^{-2}$ & $2.20\times10^{-2}$ & $3.52\times10^{-1}$ & $1.53\times10^{-1}$ \\
    \hline
    \end{tabular}}
\end{table}
\FloatBarrier

\begin{figure}[h!]
    \centering

    \begin{subfigure}[t]{0.49\textwidth}
        \centering
        \includegraphics[width=\linewidth]{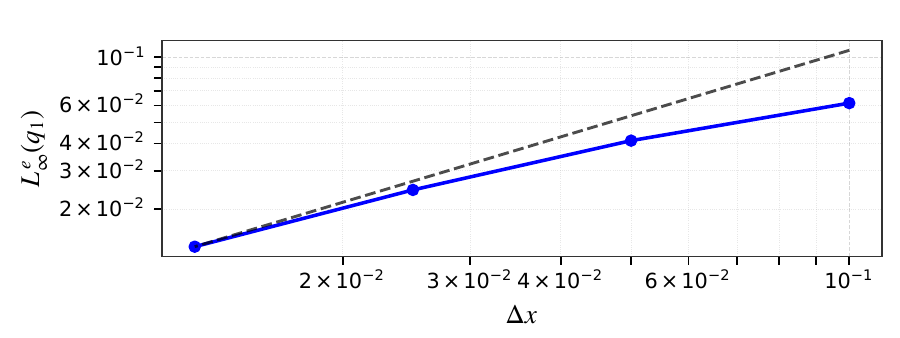}
        \caption{Virus ($q_1$)}
        \label{fig:convergence_Linf_q1}
    \end{subfigure}
    \hspace{-0.9em}
    \begin{subfigure}[t]{0.49\textwidth}
        \centering
        \includegraphics[width=\linewidth]{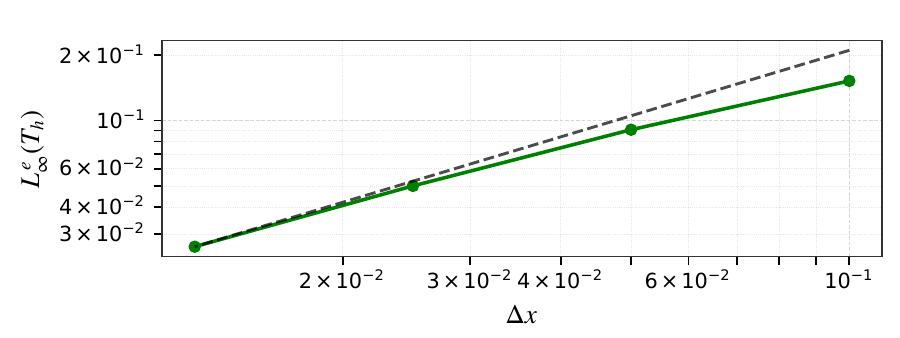}
        \caption{T helper ($T_h$)}
        \label{fig:convergence_Linf_Th}
    \end{subfigure}

    \begin{subfigure}[t]{0.49\textwidth}
        \centering
        \includegraphics[width=\linewidth]{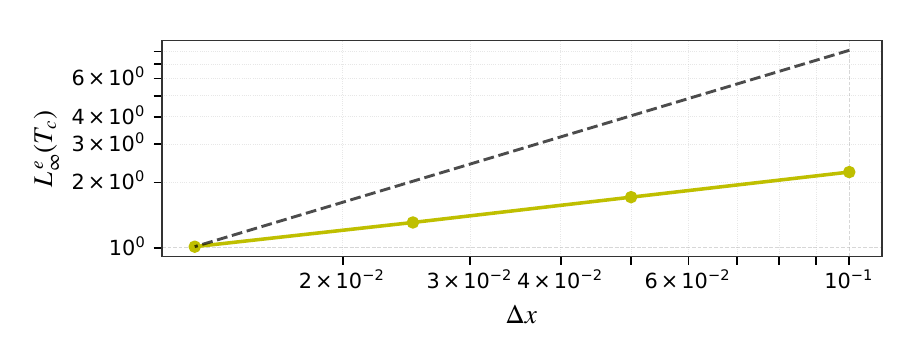}
        \caption{T cytotoxic ($T_c$)}
        \label{fig:convergence_Linf_Tc}
    \end{subfigure}
    \hspace{-0.9em}
    \begin{subfigure}[t]{0.49\textwidth}
        \centering
        \includegraphics[width=\linewidth]{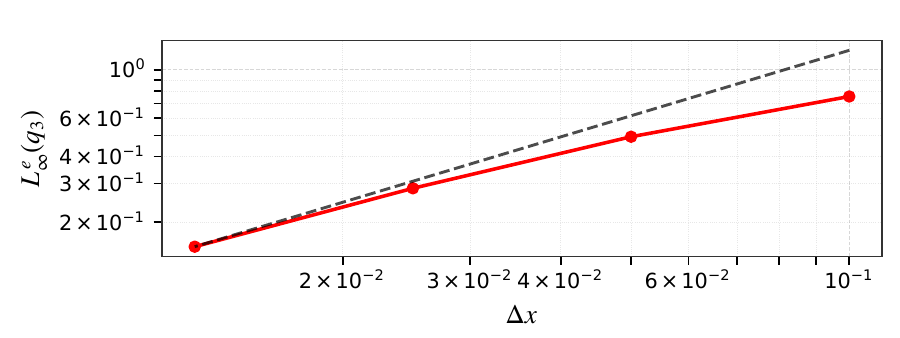}
        \caption{Cytokines ($q_3$)}
        \label{fig:convergence_Linf_q3}
    \end{subfigure}

    \vspace{0.2em}

    \includegraphics[
        width=0.48\textwidth
    ]{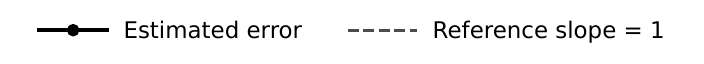}

    \vspace{-0.5em}

    \caption{Logarithmic representation of the $L_\infty^e$-norm error against $\Delta x$ for the four model variables.}
    \label{fig:validacion_convergencia_Linf_1D}
\end{figure}

\begin{figure}[h!]
    \centering

    \begin{subfigure}[t]{0.49\textwidth}
        \centering
        \includegraphics[width=\linewidth]{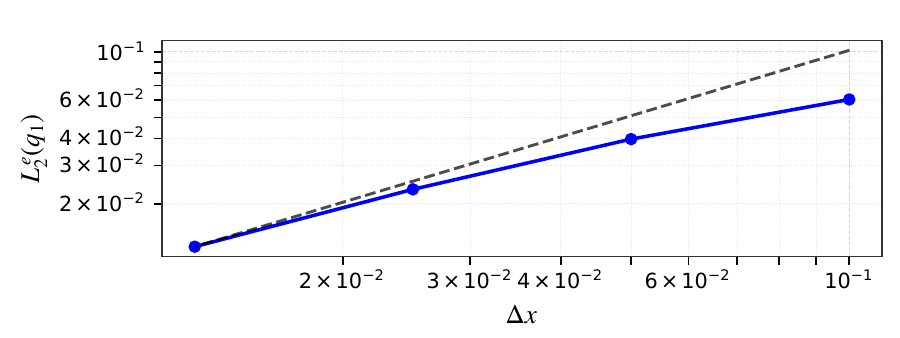}
        \caption{Virus ($q_1$)}
        \label{fig:convergence_L2_q1}
    \end{subfigure}
    \hspace{-0.9em}
    \begin{subfigure}[t]{0.49\textwidth}
        \centering
        \includegraphics[width=\linewidth]{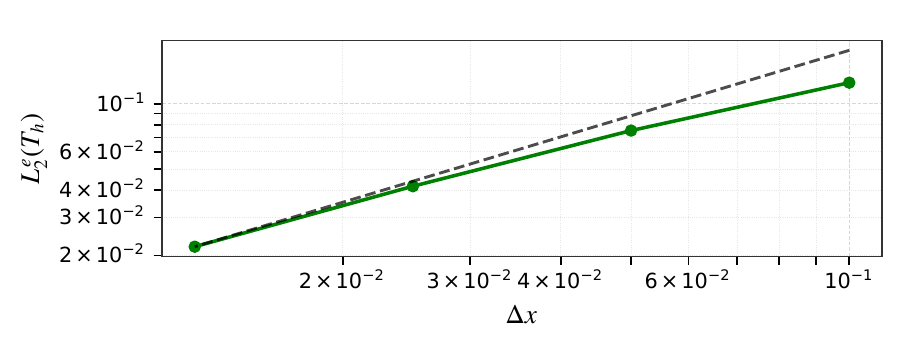}
        \caption{T helper ($T_h$)}
        \label{fig:convergence_L2_Th}
    \end{subfigure}

    \begin{subfigure}[t]{0.49\textwidth}
        \centering
        \includegraphics[width=\linewidth]{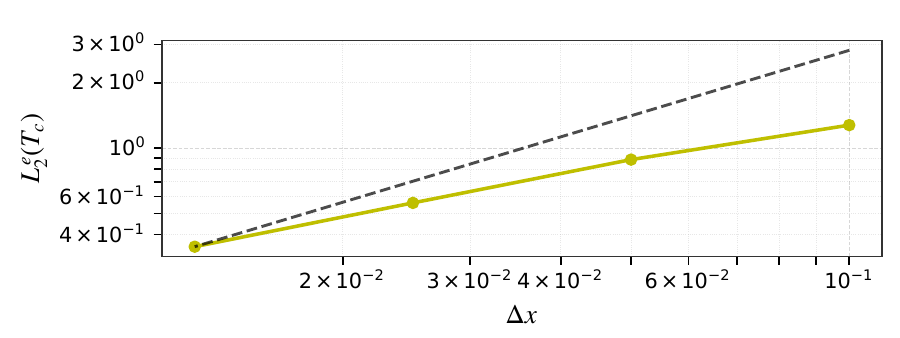}
        \caption{T cytotoxic ($T_c$)}
        \label{fig:convergence_L2_Tc}
    \end{subfigure}
    \hspace{-0.9em}
    \begin{subfigure}[t]{0.49\textwidth}
        \centering
        \includegraphics[width=\linewidth]{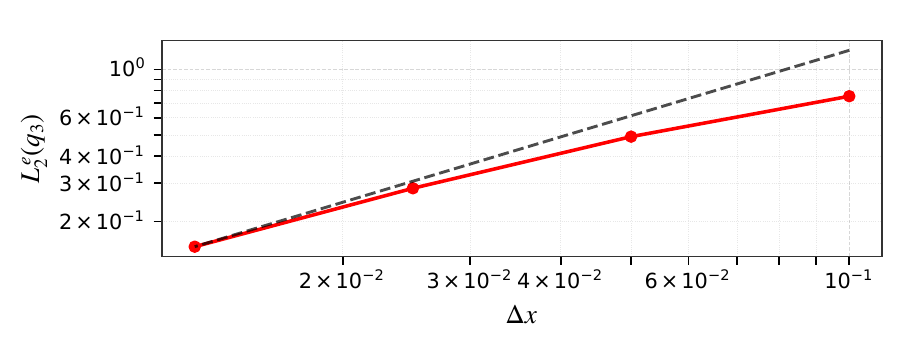}
        \caption{Cytokines ($q_3$)}
        \label{fig:convergence_L2_q3}
    \end{subfigure}

    \vspace{0.2em}

    \includegraphics[
        width=0.48\textwidth
    ]{Imagenes/CASO_1D/convergence_L2_common_legend.pdf}

    \vspace{-0.5em}

    \caption{Logarithmic representation of the $L_2^e$-norm error against $\Delta x$ for the four model variables.}
    \label{fig:validacion_convergencia_L2}
\end{figure}

The error decreases progressively with mesh refinement for all four variables and in both norms, indicating convergence. The logarithmic plots in Figures \ref{fig:validacion_convergencia_Linf_1D} and \ref{fig:validacion_convergencia_L2} show an approximately linear trend approaching first order of accuracy.
\FloatBarrier
\FloatBarrier

\begin{figure}[h!]
    \centering

    \begin{subfigure}[t]{0.49\textwidth}
        \centering
        \includegraphics[width=\linewidth]{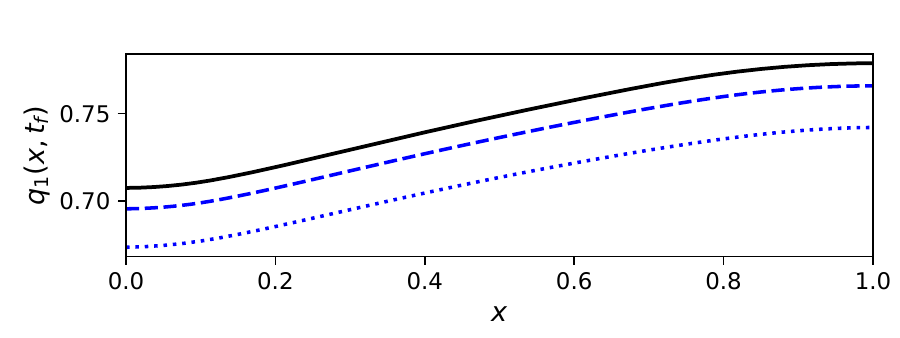}
        \caption{Virus ($q_1$)}
        \label{fig:convergence_q1}
    \end{subfigure}
    \hspace{-1.1em}
    \begin{subfigure}[t]{0.49\textwidth}
        \centering
        \includegraphics[width=\linewidth]{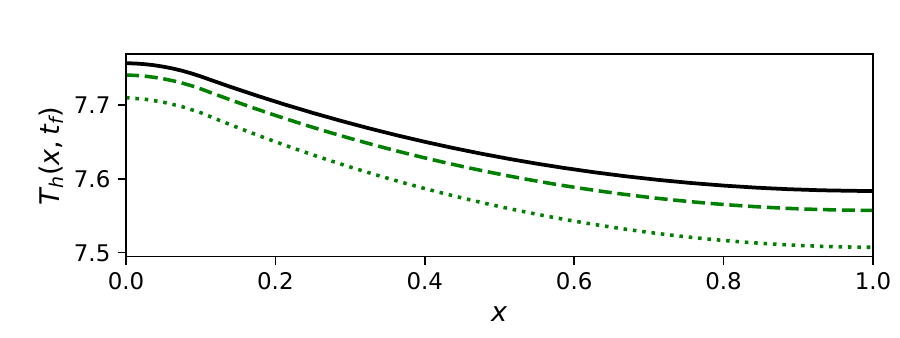}
        \caption{T helper ($T_h$)}
        \label{fig:convergence_Th}
    \end{subfigure}

    \begin{subfigure}[t]{0.49\textwidth}
        \centering
        \includegraphics[width=\linewidth]{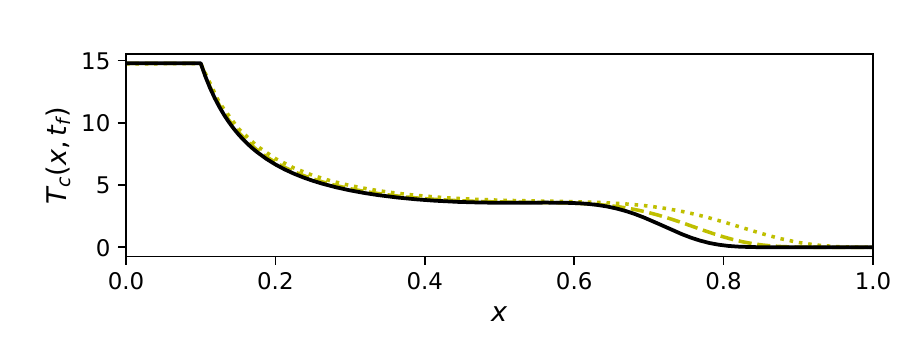}
        \caption{T cytotoxic ($T_c$)}
        \label{fig:convergence_Tc}
    \end{subfigure}
    \hspace{-1.1em}
    \begin{subfigure}[t]{0.49\textwidth}
        \centering
        \includegraphics[width=\linewidth]{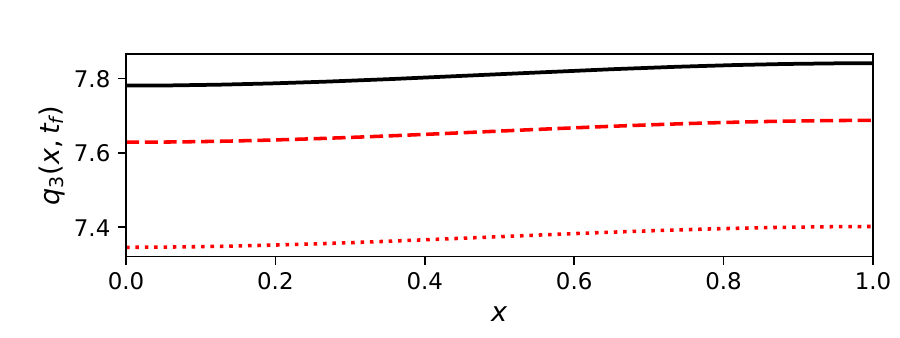}
        \caption{Cytokines ($q_3$)}
        \label{fig:convergence_q3}
    \end{subfigure}

    \vspace{0.2em}

    \includegraphics[
        width=0.48\textwidth
    ]{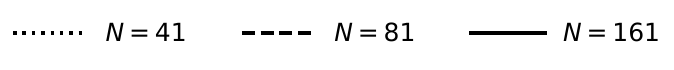}

    \vspace{-0.5em}

    \caption{Spatial mesh-convergence comparison of the four model
    variables at the final simulation time, $t_f$, using meshes with
    $N=41$, $N=81$, and $N=161$ nodes.}

    \label{fig:spatial_convergence}
\end{figure}

Figure \ref{fig:spatial_convergence} compares the spatial profiles of the four variables for $N=41$, $N=81$ and $N=161$. Results are consistent with the previous refinement study results.

\FloatBarrier

\subsection{CPU-GPU computational performance}
\label{subsec: CPU_GPU_performance}
\begin{figure}[h!]
    \centering
    \includegraphics[width=0.6\textwidth]{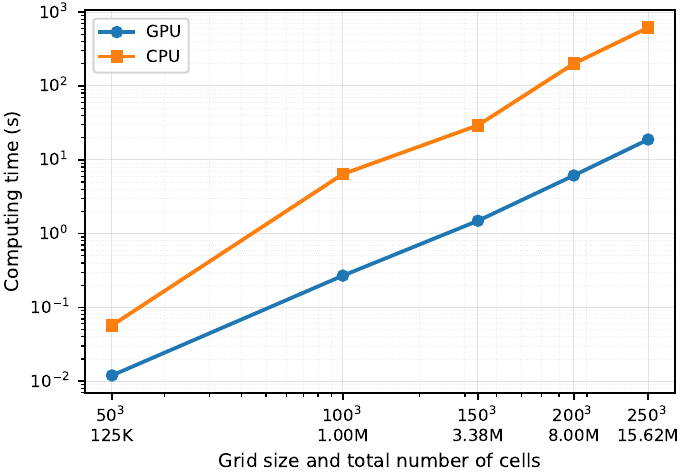}
    \caption{Comparison of CPU/GPU execution time when refining the mesh.}
    \label{fig:comparacion_cpu_gpu}
\end{figure}

In this section, CPU and GPU execution times are compared for five mesh resolutions ranging from $50 \times 50 \times 50$ to $250 \times 250 \times 250$ nodes. A detailed description of the benchmark test configuration is
provided in \ref{app:CPU_GPU_test_case}. The results are reported in Table \ref{tab:comparacion_cpu_gpu} and shown in Figure \ref{fig:comparacion_cpu_gpu}. As the mesh size increases, the speed-up of GPU versus CPU increases reaching a value around 32 for the two larger meshes. This confirms that GPU use is highly beneficial for large three-dimensional simulations, in which the number of operations increases substantially, whereas the CPU remains preferable for simulations with a lower computational load, where the computational overhead associated with GPU initialisation, data transfer, and operation launching
may be comparable to the cost of the simulation itself.
\begin{table}[h!]
\centering
\caption{Comparison of CPU and GPU execution time for different mesh resolutions.}
\label{tab:comparacion_cpu_gpu}
\renewcommand{\arraystretch}{1.1}
\makebox[\textwidth][c]{
    \resizebox{\textwidth}{!}{
    \begin{tabular}{c c c c c c}
    \hline
    \multirow{2}{*}{Mesh} &
    \multirow{2}{*}{Total cells} &
    \multirow{2}{*}{$\Delta t$} &
    \multicolumn{2}{c}{Execution time (s)} &
    \multirow{2}{*}{Speed-up} \\
    \cline{4-5}
    & & & GPU & CPU & \\
    \hline
    $50 \times 50 \times 50$    
    & $1.25 \times 10^5$ 
    & $6.94 \times 10^{-5}$   
    & 0.012  
    & 0.057
    & 4.750\\

    $100 \times 100 \times 100$ 
    & $1.00 \times 10^6$ 
    & $1.70 \times 10^{-5}$ 
    & 0.270  
    & 6.412 
    & 23.748\\

    $150 \times 150 \times 150$ 
    & $3.38 \times 10^6$ 
    & $7.51 \times 10^{-6}$  
    & 1.494  
    & 29.318 
    & 19.624\\

    $200 \times 200 \times 200$ 
    & $8.00 \times 10^6$ 
    & $4.21 \times 10^{-6}$ 
    & 6.133  
    & 200.703 
    & 32.725\\

    $250 \times 250 \times 250$ 
    & $1.56 \times 10^7$ 
    & $2.69 \times 10^{-6}$ 
    & 18.822 
    & 616.513 
    & 32.755\\
    \hline
    \end{tabular}}}
\end{table}

\section{Results}
\label{sec:results}

This section presents the analysis of the local dynamics of the model using phase space analysis, considering the regimes of chronicity and viral clearance. Subsequently, the spatial evolution of the infection is studied in 1D, 2D, and 3D configurations, progressively incorporating spatial transport and the heterogeneity of the environment.

\subsection{Phase space analysis of local dynamics}
\label{sec:phase_space_local}

As a preliminary analysis, the local dynamics of the 0D model are studied using phase space, considering only the reactions and without spatial dependence.
This analysis allows us to identify the two regimes considered throughout this work: chronicity and viral clearance. Both regimes are obtained by modifying only the parameter $a_5$, which characterizes the efficiency of cytotoxic T lymphocytes in virus clearance. The values of the parameters used are summarized in Table \ref{tab:parameters_0D}.
\begin{table}[h!]
\centering
\caption{Model parameter values for the chronic and healing regimes.}
\label{tab:parameters_0D}
\resizebox{\textwidth}{!}{
\begin{tabular}{l c c c c c c c c c c c c c c}
\hline 
& $a_1$ & $C_1$ & $\epsilon$ & $\kappa$ & $a_5$ & $a_{2h}$ & $C_{Th}$ & $a_{6h}$ & $a_{6c}$ & $C_{Tc}$ & $a_3$ & $a_{nd}$ & $a_{2c}$ \\
\hline
Chronic & 1.0 & 1.0 & 0.05 & 0.01 & \textbf{0.05} & 2.0 & 8.0 & 0.2 & 0.2 & 15.0 & 0.8 & 0.6 & 2.0 \\
Healing & 1.0 & 1.0 & 0.05 & 0.01 & \textbf{0.08} & 2.0 & 8.0 & 0.2 & 0.2 & 15.0 & 0.8 & 0.6 & 2.0 \\
\hline
\end{tabular}}
\end{table}
\FloatBarrier

In both phase portraits, the streamlines indicate the direction of the system's evolution, while the colors represent the logarithm of the velocity magnitude. This logarithmic scale allows regions of faster and slower evolution to be distinguished more clearly.

\begin{figure}[h!]
    \centering
    \includegraphics[width=0.9\textwidth]{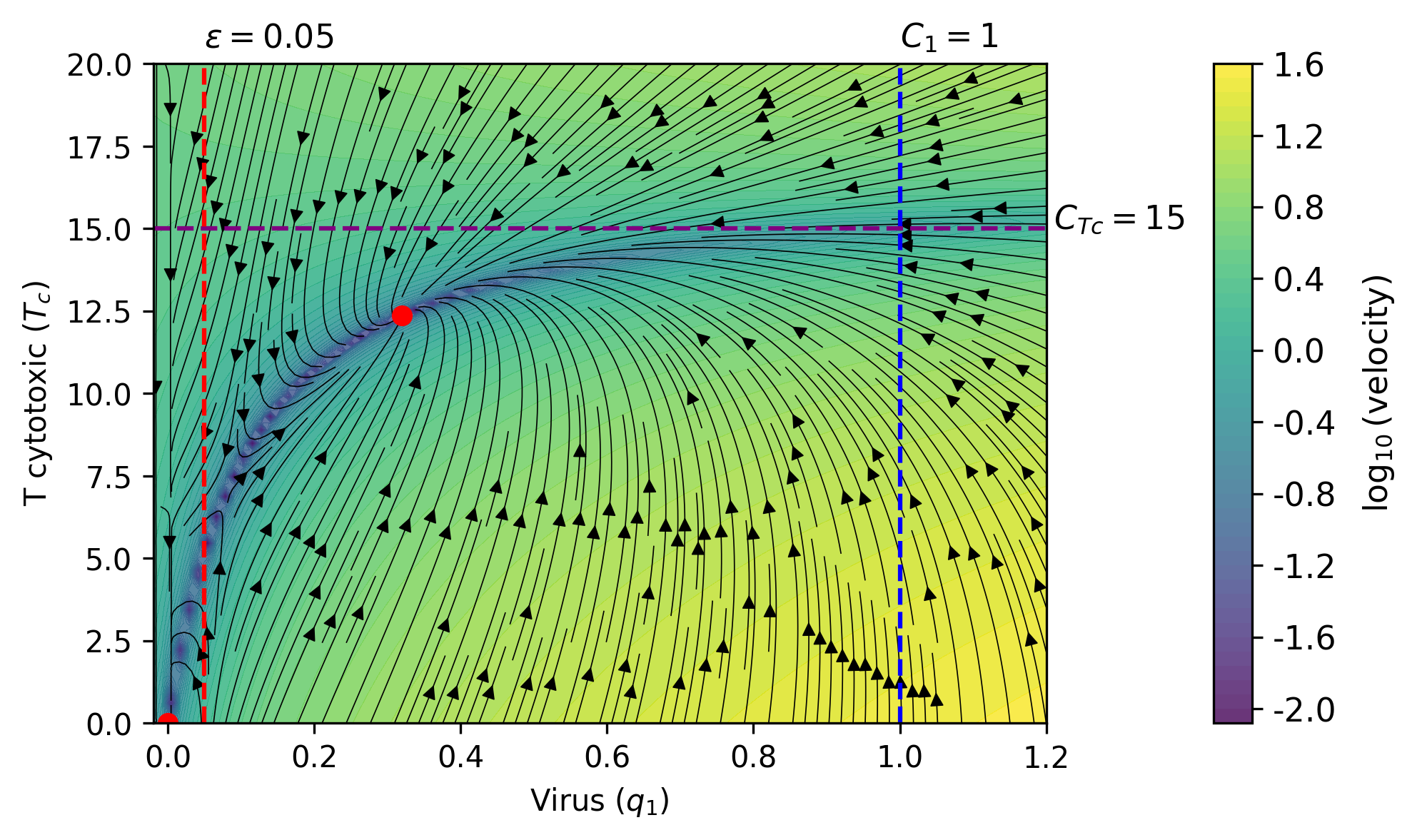}
    \caption{Velocity field for the chronic regime case, of the viral quantity (\(q_1\)) and cytotoxic T lymphocytes (\(T_c\)), the attractor being indicated by a red dot.
    The velocity represents how fast the system evolves through the phase space along the trajectories. }
    \label{fig:CasoA_campo_velocidades}
\end{figure}

In the chronic regime, the system's trajectories converge towards a stable equilibrium with a positive viral load, represented by the interior red dot in Figure \ref{fig:CasoA_campo_velocidades}. This attractor corresponds to the persistence of the infection.

\begin{figure}[h!]
    \centering
    \includegraphics[width=0.9\textwidth]{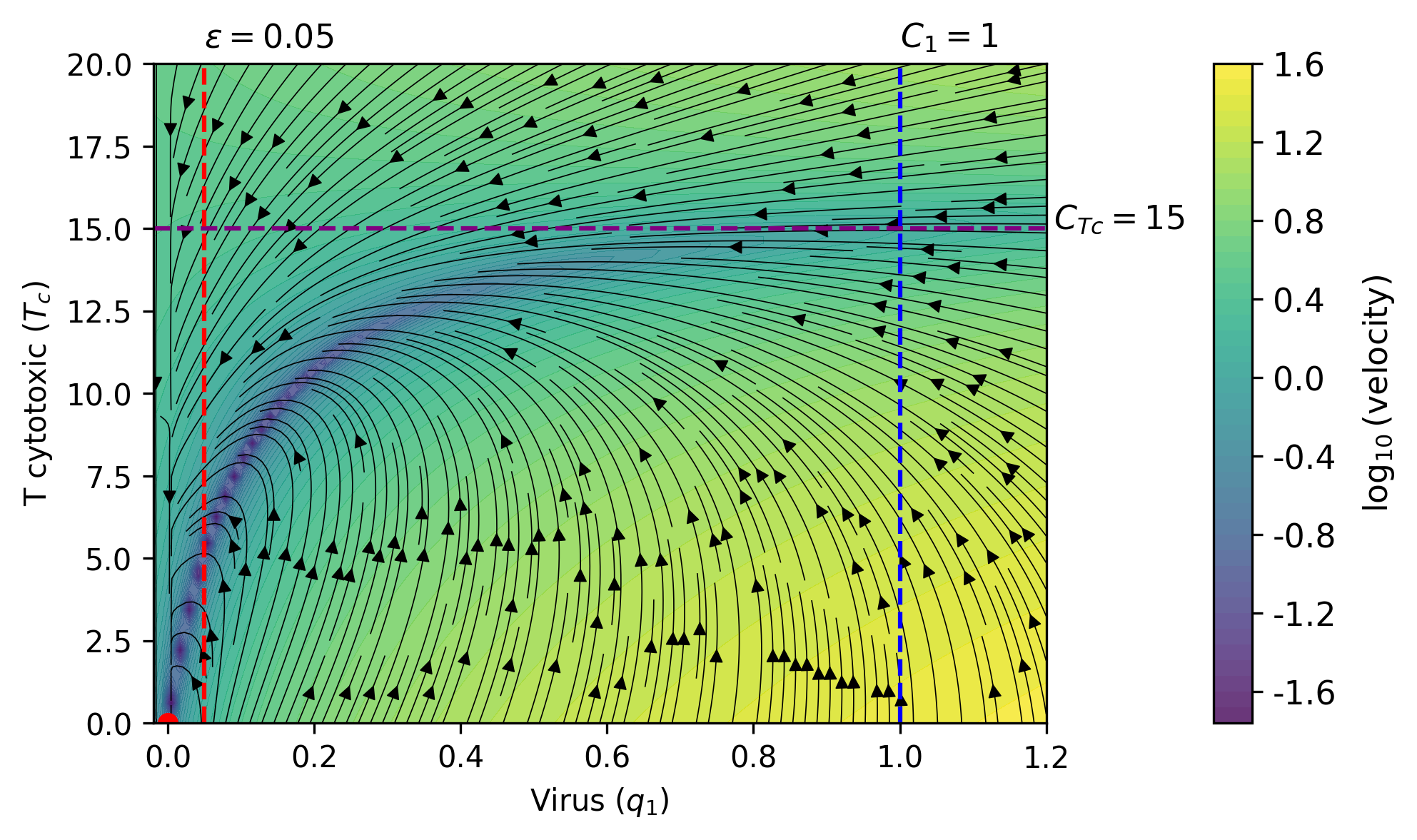}
    \caption{Velocity field for the case of the cure regime, of the viral quantity (\(q_1\)) and the cytotoxic T lymphocytes (\(T_c\)), the attractor being indicated by a red dot.}     \label{fig:CasoB_campo_velocidades}
\end{figure}

In the viral clearance regime represented in Figure \ref{fig:CasoB_campo_velocidades}, the trajectories converge towards the origin, which represents the disappearance of the viral load and, therefore, the elimination of the infection. This behavior is obtained by increasing the parameter $a_5$ with respect to the chronic case.
\FloatBarrier

\subsection{One-dimensional simulations}
The one-dimensional model was first used to analyse how the incorporation of spatial transport modifies the chronification and viral-clearance regimes identified in a previous analysis of the local reaction dynamics (0D model).
The computational domain is defined as $x \in [0, L]$, with $L=1$, and the system is integrated up to $t_f=50$ using a stability factor $\sigma=0.4$, see Equation \eqref{eq:sigma}.
Both regimes share the same initial conditions,
\begin{equation}
    q_1(x,0)=0.1, \quad T_h(x,0)=0, \quad T_c(x,0)=0, \quad q_3(x,0)=0
    \label{eq:condiciones_iniciales_1D}
\end{equation}
and the same immune-response entry region, defined over the interval $[0<x<0.1]$ through the function $\chi_{\Theta}(\mathbf{x})$. This represents a localised vascular entry point rather than a uniform source of T lymphocytes throughout the domain. The reaction parameters and transport coefficients for both regimes are reported in Tables \ref{tab:parametros_1D_reaccion} and \ref{tab:parametros_1D_dif_adv}.
\begin{table}[h!]
\centering
\caption{Model parameter values for the chronic and healing regimes.}
\label{tab:parametros_1D_reaccion}
\resizebox{\textwidth}{!}{
\begin{tabular}{l c c c c c c c c c c c c c}
\hline
 & $a_1$ & $C_1$ & $\epsilon$ & $\kappa$ & $a_5$ & $a_{2h}$ & $C_{Th}$ & $a_{6h}$ & $a_{6c}$ & $C_{Tc}$ & $a_3$ & $a_{nd}$ & $a_{2c}$ \\
\hline
Chronic   & 1.0 & 1.0 & 0.05 & 0.01 & \textbf{0.05} & 2.0 & 8.0 & 0.2 & 0.2 & 15.0 & 0.8 & 0.6 & 2.0 \\
Healing   & 1.0 & 1.0 & 0.05 & 0.01 & \textbf{0.08} & 2.0 & 8.0 & 0.2 & 0.2 & 15.0 & 0.8 & 0.6 & 2.0 \\
\hline
\end{tabular}}
\end{table} 
\FloatBarrier

\begin{table}[h!]
\centering
\caption{Model diffusion and advection coefficients for the chronic and healing regimes.}
\label{tab:parametros_1D_dif_adv}
\begin{tabular}{l c c c c}
\hline
 & $d_1^{\mathrm{ctc}}$ & $d_{T_h}^{\mathrm{diff}}$ & $d_{T_c}^{\mathrm{chem}}$ & $d_3^{\mathrm{diff}}$ \\
\hline
Chronic   & 0.6 & 0.9 & 0.5 & 0.5 \\
Healing   & 0.005 & 0.5 & 0.1 & 0.06 \\
\hline
\end{tabular}
\end{table}

\subsubsection{Chronification regime in 1D}

\begin{figure}[h!]
    \centering
    
    \begin{subfigure}[b]{0.48\textwidth}
        \centering
        \includegraphics[width=\textwidth]{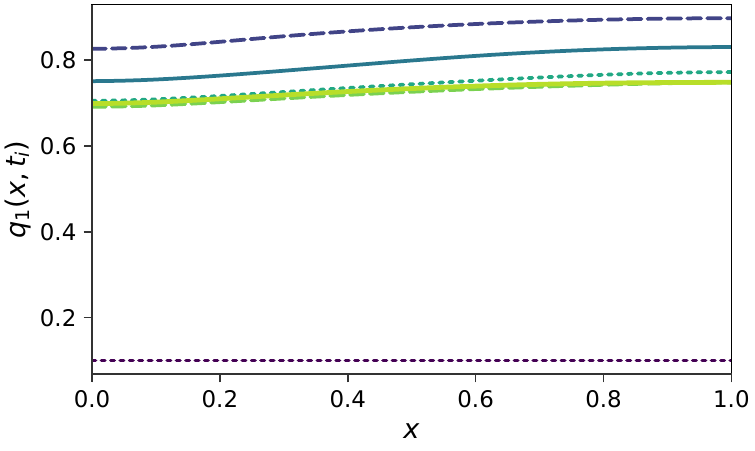}
        \caption{Virus ($q_1$).}
        \label{fig:viral_load1}
    \end{subfigure}
    \hfill
    \begin{subfigure}[b]{0.48\textwidth}
        \centering
        \includegraphics[width=\textwidth]{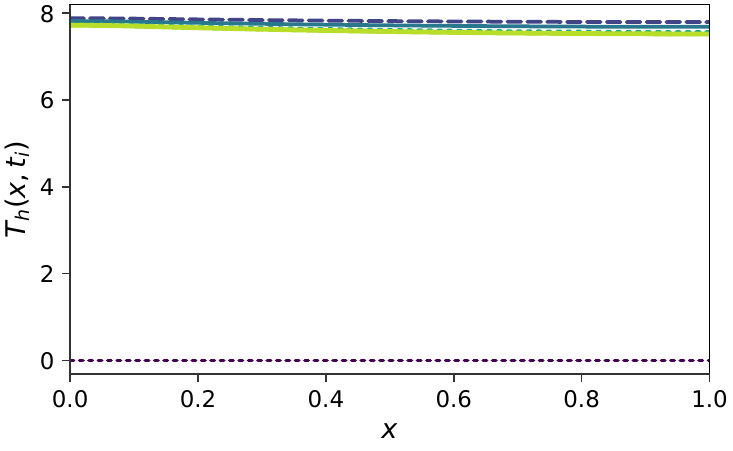}
        \caption{T helper ($T_h$).}
        \label{fig:helper_t1}
    \end{subfigure}
    
    \vspace{0.4cm} 
    \begin{subfigure}[b]{0.48\textwidth}
        \centering
        \includegraphics[width=\textwidth]{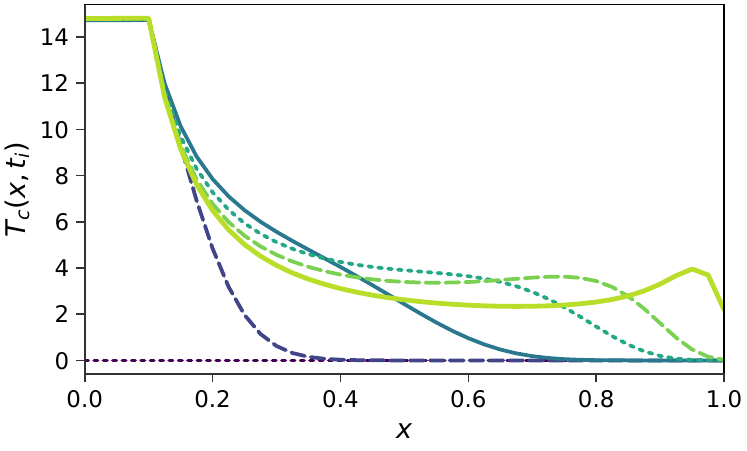}
        \caption{T cytotoxic ($T_c$).}
        \label{fig:cytotoxic_t1}
    \end{subfigure}
    \hfill
    \begin{subfigure}[b]{0.48\textwidth}
        \centering
        \includegraphics[width=\textwidth]{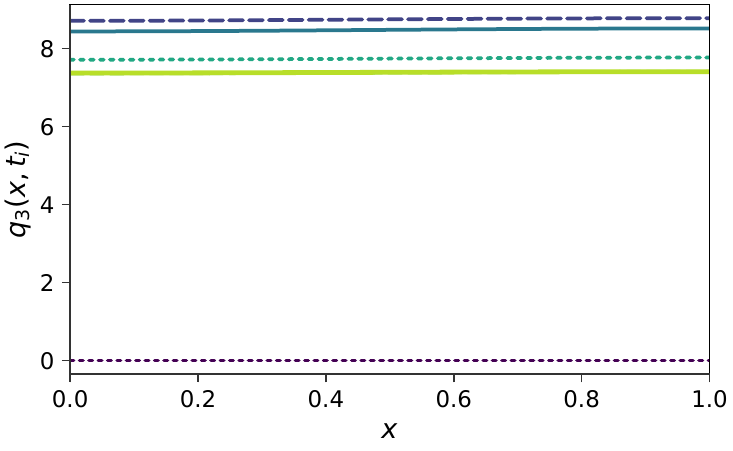}
        \caption{Cytokines $q_3$.}
        \label{fig:cytokines1}
    \end{subfigure}
    
    \vspace{0.4cm}
    
    \centering
    \includegraphics[width=0.85\textwidth]{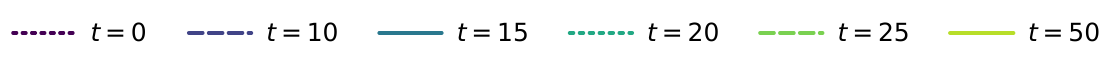}

    \caption{One-dimensional spatial plot of the chronification regime for the four model variables at different time instants $t_i$. Subfigures (a–d) share the common timeline legend shown below the panels.}
    \label{fig:CasoA_espacial_1D_0_50}
\end{figure}

Under the configuration corresponding to the chronic regime, cytotoxic T lymphocytes enter the domain through the localised inlet region $\Theta$ and propagate spatially  (Figure \ref{fig:CasoA_espacial_1D_0_50}). However, the immune response progressively weakens with distance from this region and fails to eliminate the infection throughout the domain. Consequently, the viral load evolves towards a spatially heterogeneous state with positive values that are slightly higher in regions farther from the immune-response entry point.

\begin{figure}[h!]
    \centering

    \includegraphics[
        width=\textwidth
    ]{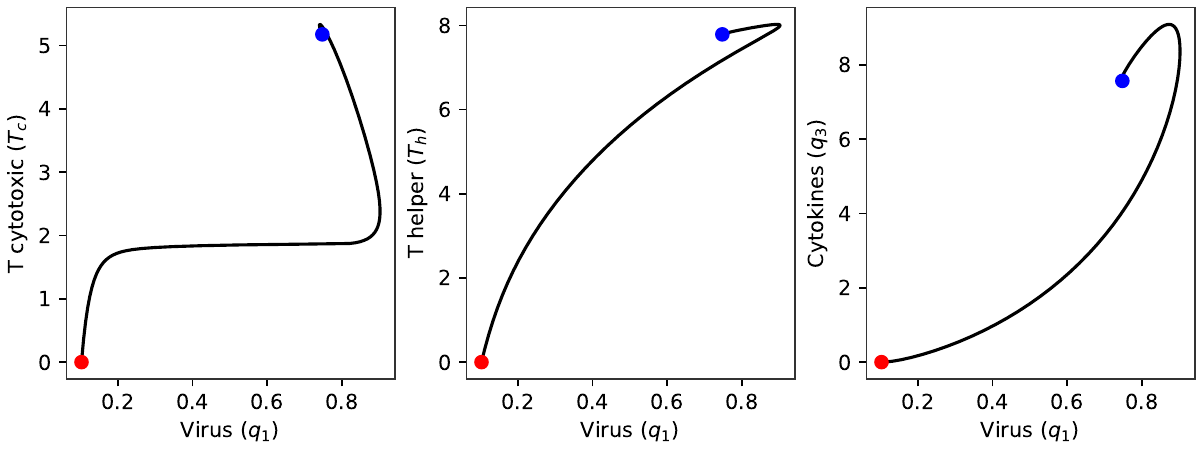}

    \vspace{0.2em}

    \includegraphics[
        width=0.45\textwidth
    ]{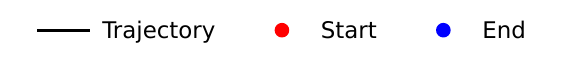}

    \vspace{-0.4em}

    \caption{Phase diagram of the $L_1(\Omega)$ norms of the immune response ($T_c$, $T_h$, and $q_3$) as a function of the $L_1(\Omega)$ norm of the viral load ($q_1$) for the chronification regime.}
    \label{fig:CasoA_fase_1D}
\end{figure}

The trajectories of viral load against the remaining variables in phase space (Figure \ref{fig:CasoA_fase_1D}) converge towards a positive equilibrium consistent with infection persistence. The phase space shows the $L_1(\Omega)$ norms of the densities, e.g. $\lVert q_1 (t) \rVert_1 = \int_\Omega \lvert q_1 (\mathbf{x},t ) \rvert \, \mathrm{d} \mathbf{x}$. 
In particular, the $T_c$--$q_1$ trajectory contains an almost horizontal segment during the transient, in which the viral load continues to increase while the $T_c$-lymphocyte response varies slowly. This behaviour reflects the lag introduced by the spatial redistribution of the variables before the steady state is reached.
\FloatBarrier

\subsubsection{Healing regime in 1D}
In contrast, the healing regime differs from the chronic case in the cytotoxic elimination rate $a_5$, which increases from 0.05 to 0.08 (Table \ref{tab:parametros_1D_reaccion}), and in the diffusion and advection coefficients listed in Table \ref{tab:parametros_1D_dif_adv}.
Notably, the viral diffusion coefficient $d_1^{\mathrm{ctc}}$ is reduced by almost two orders of magnitude, limiting viral propagation throughout the domain. This process may be seen as a model for antiviral medications slowing down the cell-to-cell transmission of the virus. 

\begin{figure}[h!]
    \centering
    
    \begin{subfigure}[b]{0.48\textwidth}
        \centering
        \includegraphics[width=\textwidth]{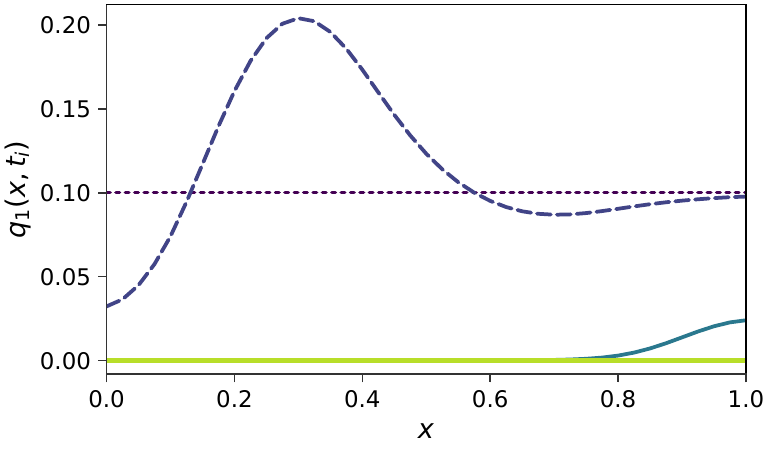}
        \caption{Virus ($q_1$).}
        \label{fig:viral_load}
    \end{subfigure}
    \hfill
    \begin{subfigure}[b]{0.48\textwidth}
        \centering
        \includegraphics[width=\textwidth]{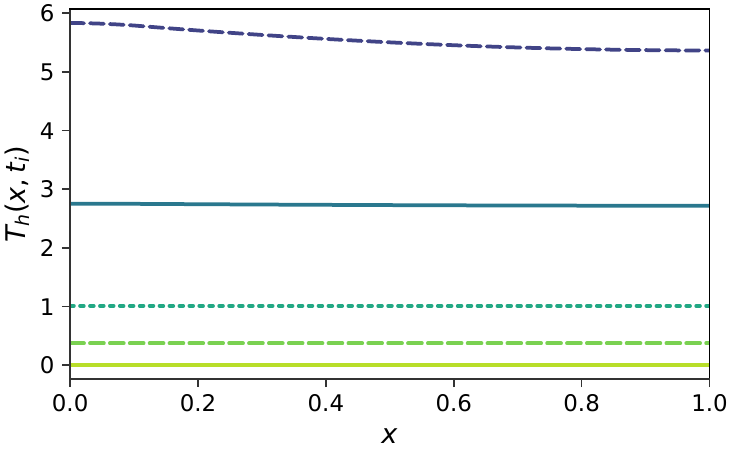}
        \caption{T helper ($T_h$).}
        \label{fig:helper_t}
    \end{subfigure}
    
    \vspace{0.4cm} 
    \begin{subfigure}[b]{0.48\textwidth}
        \centering
        \includegraphics[width=\textwidth]{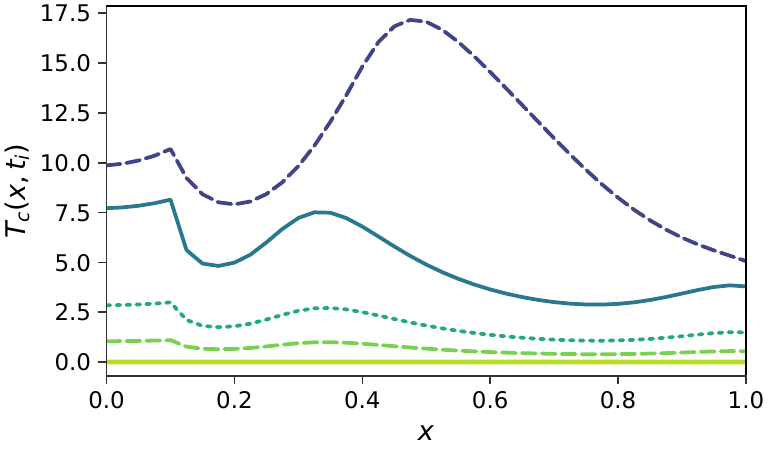}
        \caption{T cytotoxic  ($T_c$).}
        \label{fig:cytotoxic_t}
    \end{subfigure}
    \hfill
    \begin{subfigure}[b]{0.48\textwidth}
        \centering
        \includegraphics[width=\textwidth]{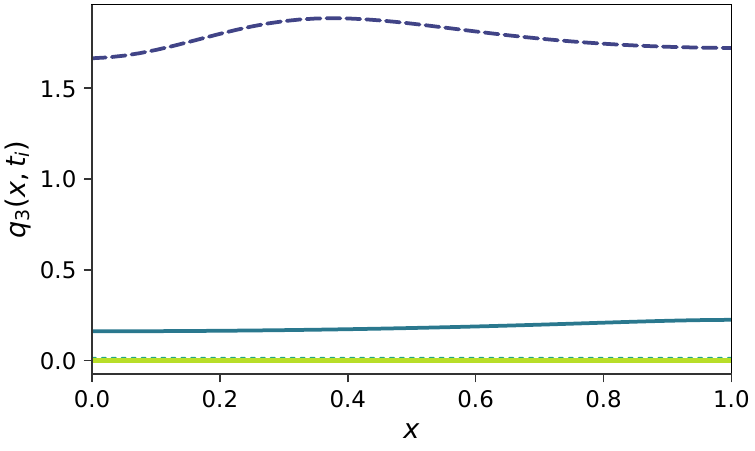}
        \caption{Cytokines ($q_3$).}
        \label{fig:cytokines}
    \end{subfigure}
    
    \vspace{0.4cm}
    
    \centering
    \includegraphics[width=0.85\textwidth]{Imagenes/CASO_1D/leyenda_comun.pdf}

    \caption{One-dimensional spatial plot of the elimination regime for the four model variables at different time instants $t_i$. Subfigures (a–d) share the common timeline legend shown below the panels.}
    \label{fig:CasoB_espacial_1D_0_50}
\end{figure}

Under this configuration, the cytotoxic T-lymphocyte wave progressively reaches the entire domain (Figure \ref{fig:CasoB_espacial_1D_0_50}), enabling elimination of the viral load. Once the infection has been cleared, the variables associated with the immune response also decrease to virtually zero by the end of the simulation.
\FloatBarrier

\begin{figure}[h!]
    \centering

    \includegraphics[
        width=\textwidth
    ]{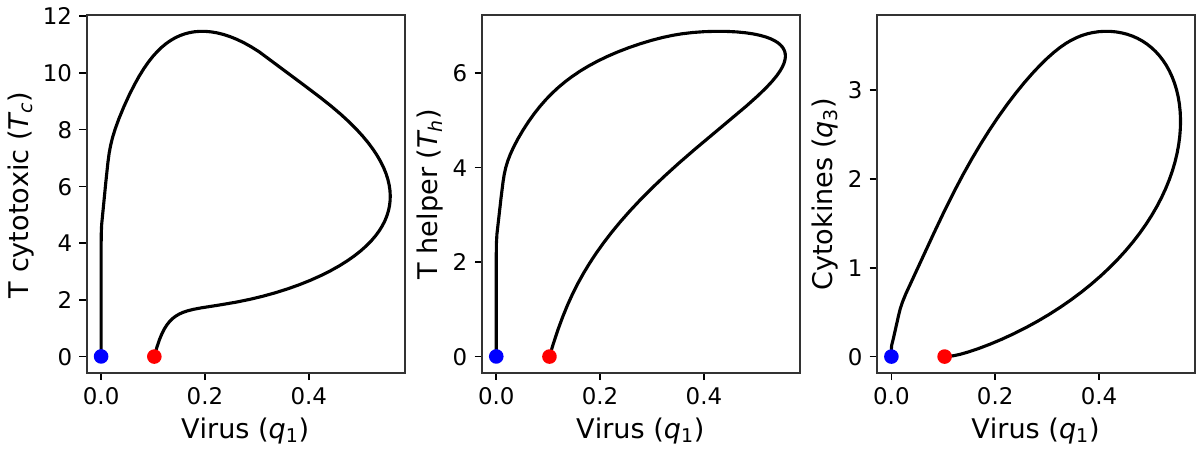}

    \vspace{0.2em}

    \includegraphics[
        width=0.45\textwidth
    ]{Imagenes/CASO_1D/leyenda_diagrama_fase.pdf}

    \vspace{-0.4em}

    \caption{Phase portraits in $L_1(\Omega)$ of the one-dimensional healing regime:
    virus against cytotoxic T cells, helper T cells, and cytokines.
    The black curves represent the temporal trajectories, while the
    red and blue markers indicate the initial and final states,
    respectively.}

    \label{fig:phase_portrait_1D_elimination}
\end{figure}

The phase space (Figure \ref{fig:phase_portrait_1D_elimination}) shows the trajectories converging towards the origin, which corresponds to infection clearance.
Unlike the chronic regime, the $T_c$--$q_1$ trajectory does not contain the almost horizontal segment observed previously, indicating a more closely coupled evolution of the cytotoxic response and the decrease in viral load.

\subsection{Two-dimensional simulations}

\subsubsection{Influence of spatial heterogeneity}\label{sec:barrier}

\begin{figure}[h!]
    \centering
    \includegraphics[width=0.6\textwidth]{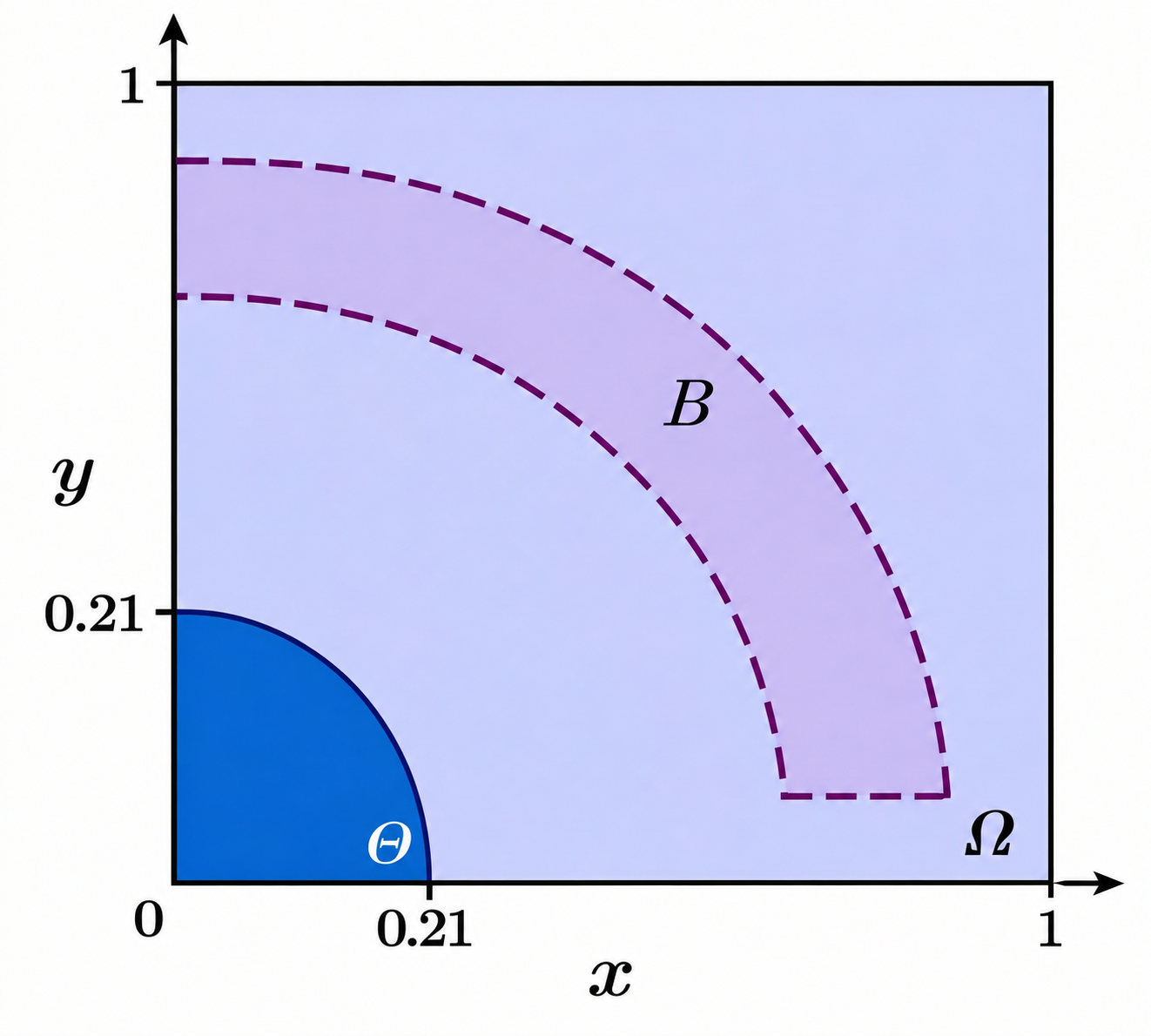}
    \caption{Schematic representation of the two-dimensional computational domain $\Omega$, showing the localized entry domain \(   \Theta\) of immune cells 
    and the internal barrier $B$.}
    \label{fig:esquema_configuracion_con_y_sin}
\end{figure}

The influence of spatial heterogeneity is analysed using two two-dimensional configurations: a homogeneous domain without internal obstacles and the same domain containing a curved barrier $B$ with a small opening (Figure \ref{fig:esquema_configuracion_con_y_sin}). The second configuration introduces a localised transport restriction without completely isolating the two regions of the domain, allowing analysis of how a spatial obstruction modifies infection propagation and, in particular, immune-response access to the regions behind the barrier.
In the medical context, these cases correspond to studying the influence of scar tissue or even fibrosis in the liver structure. Further, it is an experimental environment for testing the influence of the liver lobes on the spread of the virus. 

Heterogeneity is introduced by setting the diffusion and advection coefficients, as well as the model variables, to zero within the barrier, so that this region acts as an obstacle to spatial transport. The opening is therefore the only direct connection between the two sides of the obstacle.

The simulations were performed on a domain $x,y\in[0,L]$, with $L=L_x=L_y=1$, up to $t_f=400$ using a stability factor $\sigma=0.15$. Although the problem is analysed as two-dimensional, the implementation retains the library's general structure and uses a minimal dimension in the $z$ direction, discretised with three nodes. The same initial conditions are used in both configurations,
\begin{equation}
    q_1(\mathbf{x},0)=0.1, \quad T_h(\mathbf{x},0)=0, \quad T_c(\mathbf{x},0)=0, \quad q_3(\mathbf{x},0)=0.
    \label{eq:condiciones_iniciales_2D}
\end{equation}

The immune response is introduced through a localised region $\Theta$ at the origin, whose cross-section in the $xy$-plane is given by
$R_{xy}\leq R_{\mathrm{interior}}$, with
$R_{\mathrm{interior}}=0.21$.
In this region, $T_h$ and $T_c$ lymphocytes enter. 
The barrier $B$ is defined as an annular sector between
$R_{\mathrm{barrier,in}}=0.70$ and
$R_{\mathrm{barrier,out}}=0.90$.
The opening spans the angular interval $0\leq\theta<0.075\,\mathrm{rad}$,
measured from the positive $x$-axis, and extends across the full radial
thickness of the barrier, allowing both the virus and the variables
associated with the immune response to pass through the obstacle.

The reaction parameters and diffusion and advection coefficients are identical in both configurations and are reported in Tables \ref{tab:parametros_reaccion_2D} and \ref{tab:parametros_difusion_adveccion_2D}. The barrier is therefore the only difference between the two simulations, allowing the effect of geometric heterogeneity on the space--time evolution of the system to be isolated.
\begin{table}[h!]
\centering
\caption{Reaction parameter values for the two-dimensional analysis, common to both configurations (with and without a barrier).}
\label{tab:parametros_reaccion_2D}
\resizebox{\textwidth}{!}{
\begin{tabular}{l c c c c c c c c c c c c c}
\hline
 & $a_1$ & $C_1$ & $\epsilon$ & $\kappa$ & $a_5$ & $a_{2h}$ & $C_{Th}$ & $a_{6h}$ & $a_{6c}$ & $C_{Tc}$ & $a_3$ & $a_{nd}$ & $a_{2c}$ \\
\hline
Value   & 1.0 & 1.0 & 0.05 & 0.05 & 0.08 & 25.0 & 8.0 & 0.2 & 0.02 & 15.0 & 0.6 & 1.0 & 45.0 \\
\hline
\end{tabular}}
\end{table}
\FloatBarrier
\begin{table}[h!]
\centering
\caption{Diffusion and advection coefficient values for the two-dimensional analysis, common to both configurations (with and without a barrier).}
\label{tab:parametros_difusion_adveccion_2D}
\begin{tabular}{l c c c c}
\hline
 & $d_1^{\mathrm{ctc}}$ & $d_{T_h}^{\mathrm{diff}}$ & $d_{T_c}^{\mathrm{chem}}$ & $d_3^{\mathrm{diff}}$ \\
\hline
Value  & 0.01 & 0.5 & 0.1 & 0.06 \\
\hline
\end{tabular}
\end{table}
\FloatBarrier

In the absence of internal barriers or obstacles, the immune response can propagate throughout the domain without geometric restrictions. With no physical restrictions, the infection spreads freely throughout the domain and is ultimately eliminated by the immune response.

In particular, starting in the inlet region $ \Theta$ the immune response propagates radially, driven by diffusion and the effect of the advective term, until it covers almost the entire domain at $t=10$. From that time onwards, the immune response acts globally on the active viral load, which stops growing and decreases uniformly until it disappears completely at the final simulation time $t=400$.

\begin{figure}[h!]
    \centering
    \includegraphics[width=1.0\textwidth]{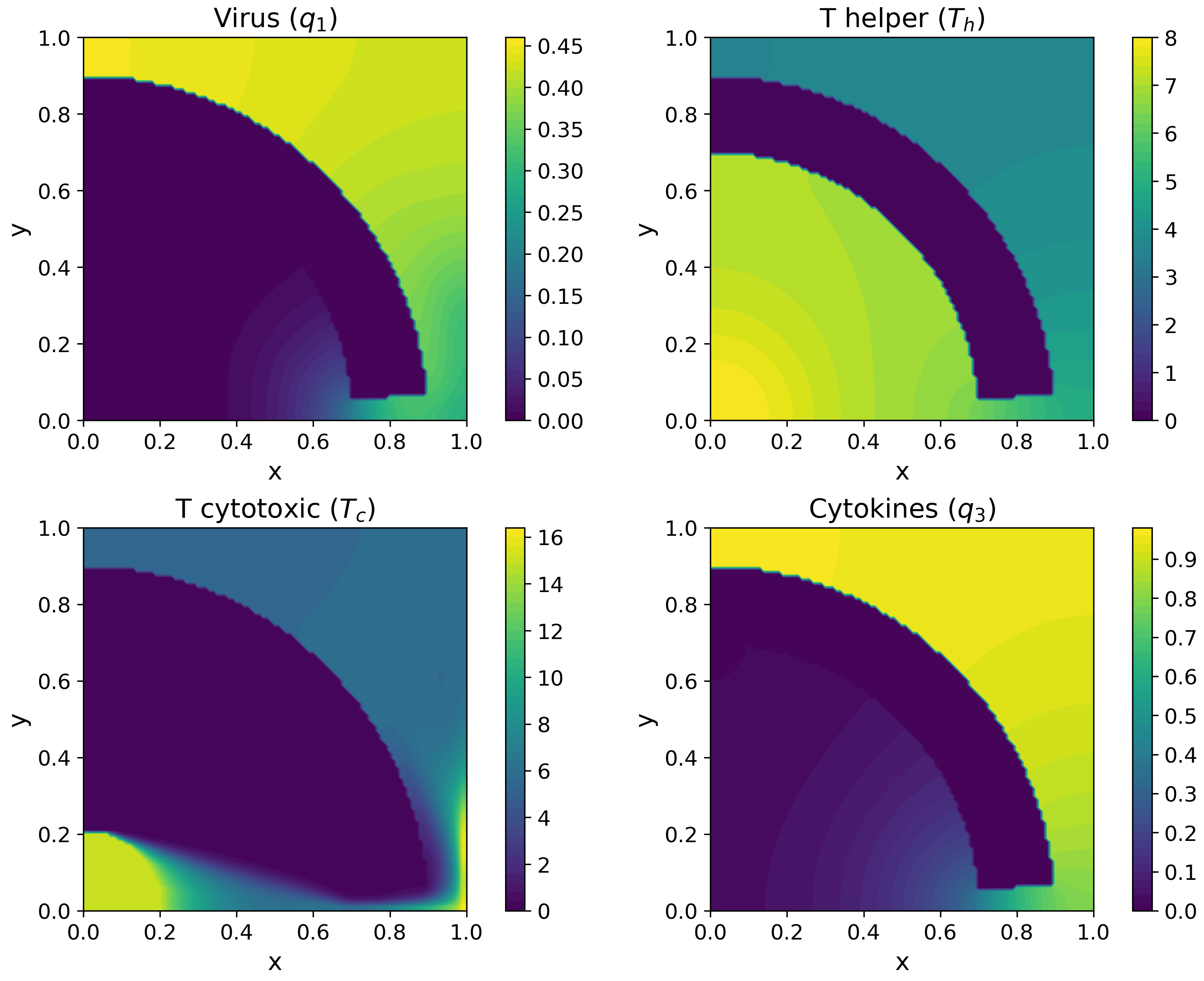}
    \caption{Spatial representation of the system dynamics in the $xy$ domain at the final time $t_f=400$ after introduction of the barrier.}
    \label{fig:representacion_espacial_2d_modelo1_barrera_400s}
\end{figure}

Introducing the barrier $B$ leads to a different spatial distribution and final state. At $t_f=400$, the viral load has been eliminated in the region before the obstacle but remains positive in part of the region behind the barrier (Figure \ref{fig:representacion_espacial_2d_modelo1_barrera_400s}). The immune response also remains active in these regions, giving rise to a spatially heterogeneous persistent state.

Initially, the $T_c$ front propagates from the inlet region $\Theta$ and reduces the viral load in the area before the barrier. As it approaches the obstacle, its propagation is restricted and the intensity of the cytotoxic response decreases. Although part of the immune response passes through the opening, it does not reach the entire domain with sufficient intensity, allowing the viral load to persist in the region behind the barrier.

\begin{figure}[htbp]
    \centering

    \begin{minipage}[b]{0.48\textwidth}
        \centering
        \includegraphics[width=\textwidth]
        {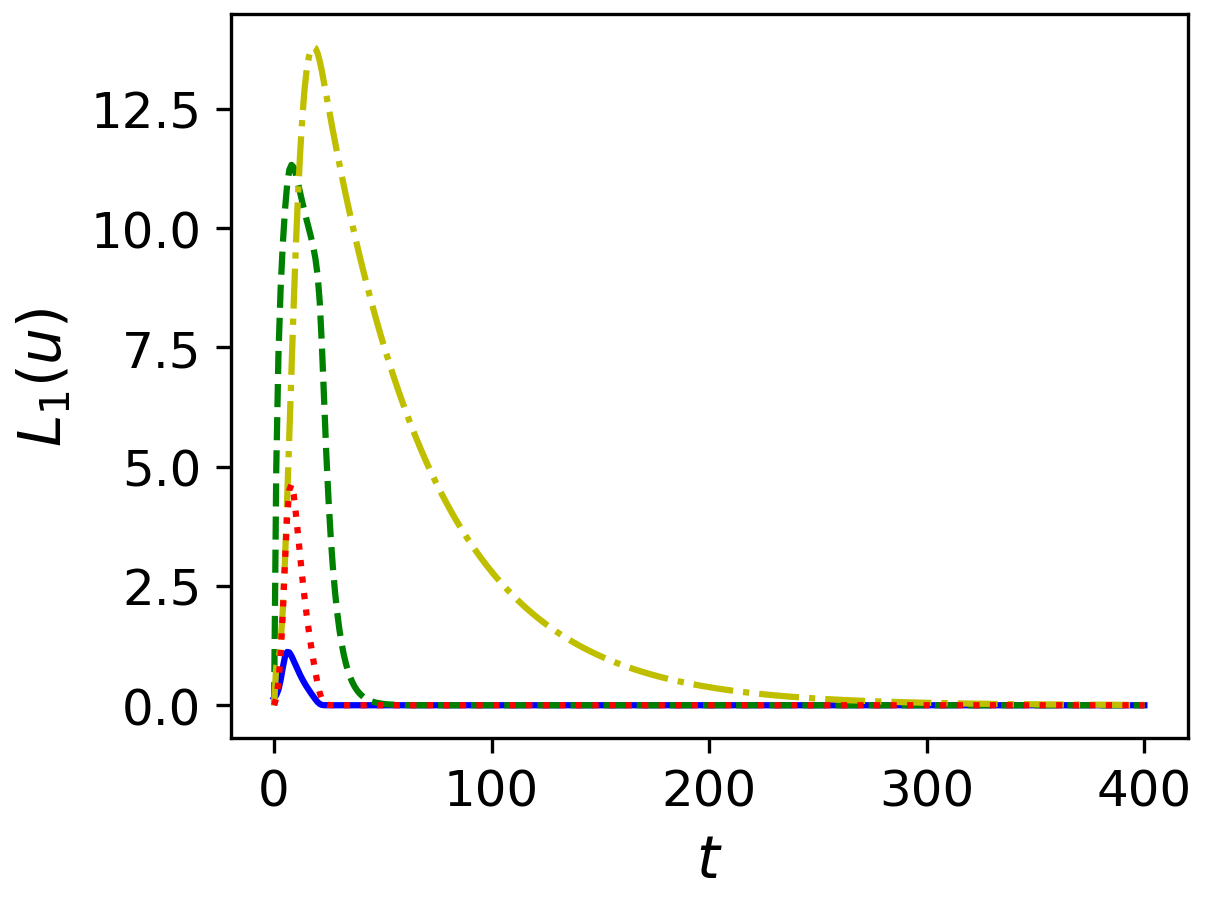}
        \par\smallskip
        (a) Without a barrier
    \end{minipage}
    \hfill
    \begin{minipage}[b]{0.48\textwidth}
        \centering
        \includegraphics[width=\textwidth]
        {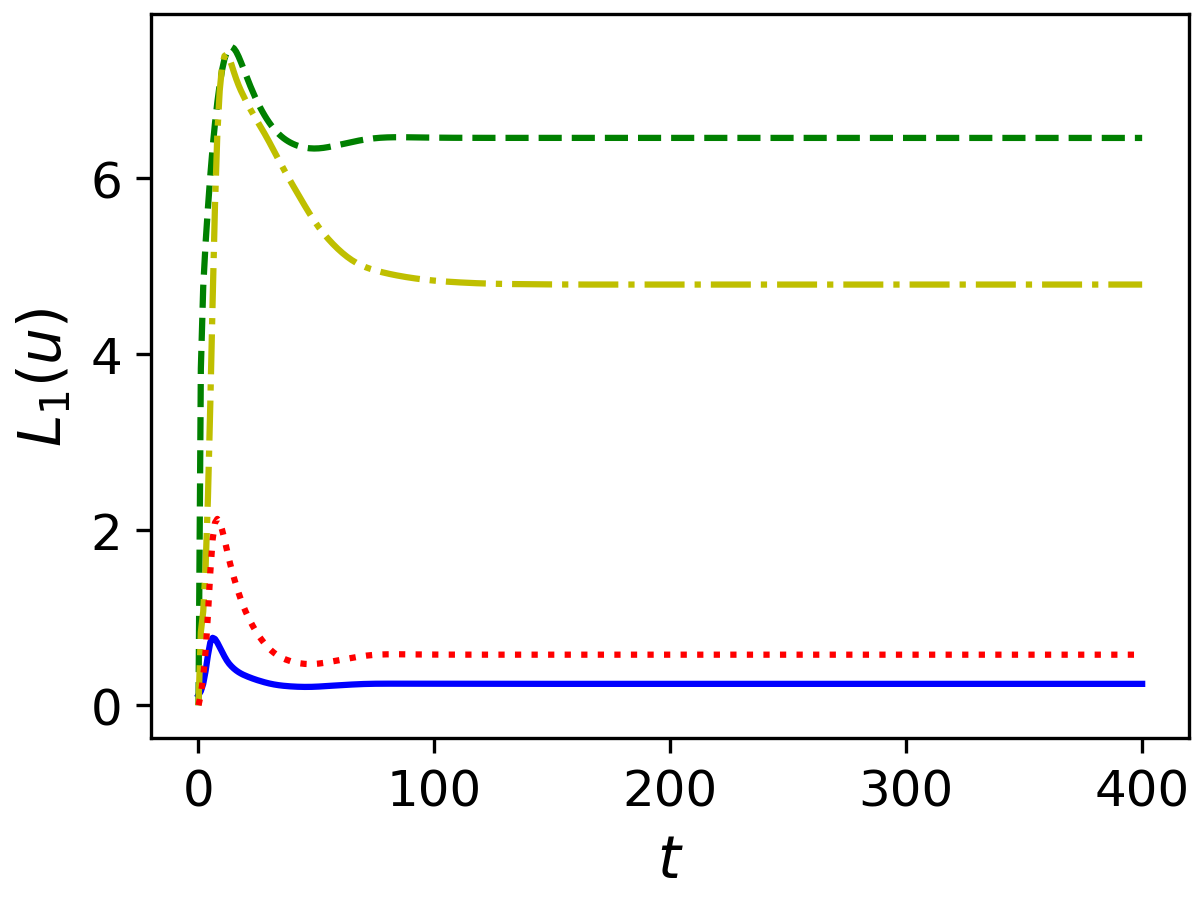}
        \par\smallskip
        (b) With a barrier
    \end{minipage}

    \par\medskip
    \includegraphics[width=0.48\textwidth]
    {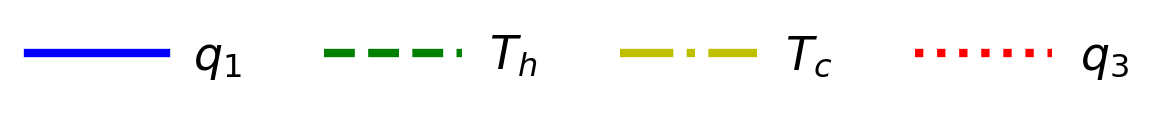}

    \caption{Temporal evolution of the spatial $L_1$ norm of the four model variables for the two-dimensional configurations without and with an internal barrier.}
    \label{fig:comparacion_barrera}
\end{figure}

Comparison of the temporal evolution of the spatially integrated densities (Figure \ref{fig:comparacion_barrera}) confirms the different regimes reached. Without a barrier, the viral load tends to zero and the immune response progressively disappears. In the presence of the barrier, by contrast, the viral load converges towards a positive state, $\lVert q_1^* \rVert_1 (\Omega) \approx0.25$, accompanied by a persistent immune response. Because both simulations share the same initial conditions, reaction parameters, and transport coefficients, the comparison shows that, within the model, introducing the barrier is sufficient to change the long-term regime from viral clearance to persistent infection.
\FloatBarrier

\subsubsection{Influence of chemotaxis in heterogeneous settings}
To isolate the role of the cytokine $q_3$ gradient in $T_c$-lymphocyte transport, a two-dimensional maze-shaped geometry was considered. The immune-response inlet and the infected region were placed at opposite ends of the domain, while the maze walls ($\Omega_{\mathrm{in}}$) were defined as impermeable regions in which both the model variables and the diffusion and advection coefficients were set to zero.
The initial conditions were defined independently inside the maze paths ($\Omega_{\mathrm{in}}$) and walls ($\Omega_{\mathrm{out}}$) as follows:
\begin{align*}
q_1(\mathbf{x},0)&=0, &
T_h(\mathbf{x},0)&=1, &
T_c(\mathbf{x},0)&=0.5, &
q_3(\mathbf{x},0)&=0.05,
&& \mathbf{x}\in\Omega_{\mathrm{in}},\\
q_1(\mathbf{x},0)&=1, &
T_h(\mathbf{x},0)&=0, &
T_c(\mathbf{x},0)&=0, &
q_3(\mathbf{x},0)&=0,
&& \mathbf{x}\in\Omega_{\mathrm{out}},
\end{align*}
with all four variables initially zero throughout the remainder of the domain. The reaction parameters and transport coefficients used are reported in Tables \ref{tab:parametros_laberinto} and \ref{tab:transporte_laberinto}.
\begin{table}[h!]
\centering
\caption{Reaction parameters used in the maze experiment.}
\label{tab:parametros_laberinto}
\resizebox{\textwidth}{!}{
\begin{tabular}{l c c c c c c c c c c c c c}
\hline
& $a_1$ & $C_1$ & $\epsilon$ & $\kappa$ & $a_5$ & $a_{2h}$ & $C_{Th}$ &
$a_{6h}$ & $a_{6c}$ & $C_{Tc}$ & $a_3$ & $a_{nd}$ & $a_{2c}$ \\
\hline
Value & 1.0 & 1.0 & 0.05 & 0.05 & 0.01 & 50.0 & 8.0 & 0.2 & 0.002 & 15.0 & 5.0 & 0.005 & 80.0 \\
\hline
\end{tabular}}
\end{table}
\FloatBarrier
\begin{table}[h!]
\centering
\caption{Transport coefficients used in the maze experiment.}
\label{tab:transporte_laberinto}
\begin{tabular}{l c c c c}
\hline
& $d_1^{\mathrm{ctc}}$ & $d_{T_h}^{\mathrm{diff}}$ &
$d_{T_c}^{\mathrm{chem}}$ & $d_3^{\mathrm{diff}}$ \\
\hline
Value & 0.0 & 0.08 & 0.15 & 0.06 \\
\hline
\end{tabular}
\end{table}
\FloatBarrier

In particular, $d_1^{\mathrm{ctc}}=0$ is considered, so the virus cannot move through the domain. The infection therefore remains localised in its initial region, allowing a more isolated analysis of the chemotactic transport of $T_c$ induced by the $q_3$ gradient. 

\begin{figure}[h!]
    \centering

    \begin{subfigure}{0.95\textwidth}
        \centering
        \includegraphics[width=\textwidth]
        {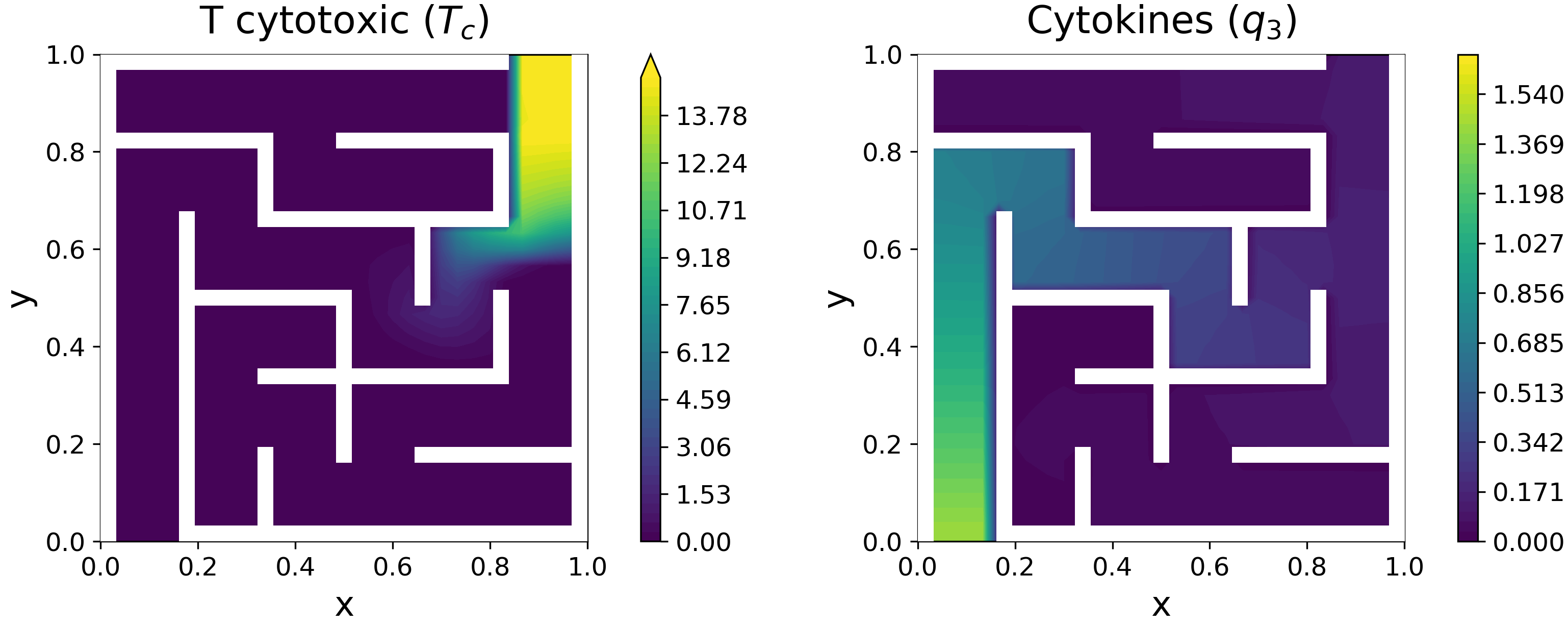}
        \caption{$t=16$}
        \label{fig:laberinto_t16}
    \end{subfigure}

    \vspace{0.3cm}

    \begin{subfigure}{0.95\textwidth}
        \centering
        \includegraphics[width=\textwidth]
        {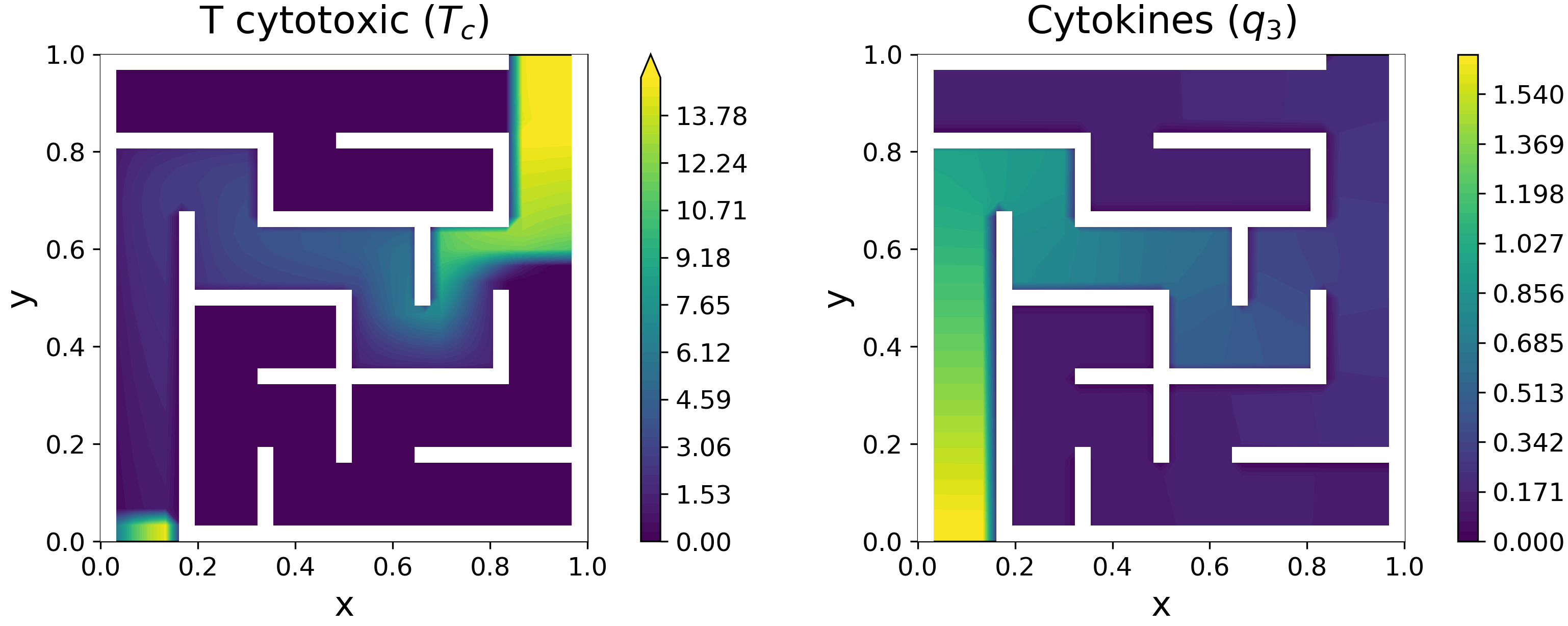}
        \caption{$t=23$}
        \label{fig:laberinto_t23}
    \end{subfigure}

    \caption{Evolution of cytotoxic T cells ($T_c$) and cytokines ($q_3$)
    in the maze at two representative time instants.}
    \label{fig:laberinto_simulacion}

\end{figure}

Figure \ref{fig:laberinto_simulacion} shows the evolution of $T_c$ and $q_3$ at the intermediate times $t=\{10,16,20,23\}$. 
Cytokines progressively diffuse through the accessible regions of the domain, generating a spatial gradient that directs the movement of cytotoxic T lymphocytes. Initially concentrated in the inlet region, the $T_c$ cells advance along this gradient and adapt to the maze geometry until they reach the infected region at approximately $t=23$, after which the viral load begins to decrease. These results show that the chemotactic mechanism directs $T_c$ transport through spatially complex geometries.

\section{Application cases}
\label{sec:3dappl}

Based on the observations for heterogeneous domains and the influence of the parameters on the disease progression, we now present simulations in a realistic three-dimensional domain. 

\subsection{Construction of the three-dimensional liver geometry}
\label{subsec:liver_geometry}

The three-dimensional computational domain was constructed from computed tomography (CT) images using 3D Slicer \cite{Fedorov2012}. Specifically, an anony\-mized whole-body CT volume was used from the public KUNAS (Kaunas University of Technology) clinical dataset, which is part of the CT-ORG multi-organ segmentation benchmark hosted on The Cancer Imaging Archive (TCIA) \cite{Rister2020}. The objective was to obtain a geometry that preserved the anatomical structure of the liver while allowing for the assignment of spatially heterogeneous transport properties to different regions.

\begin{figure}[h!]
    \centering
    \includegraphics[width=1.0\textwidth]{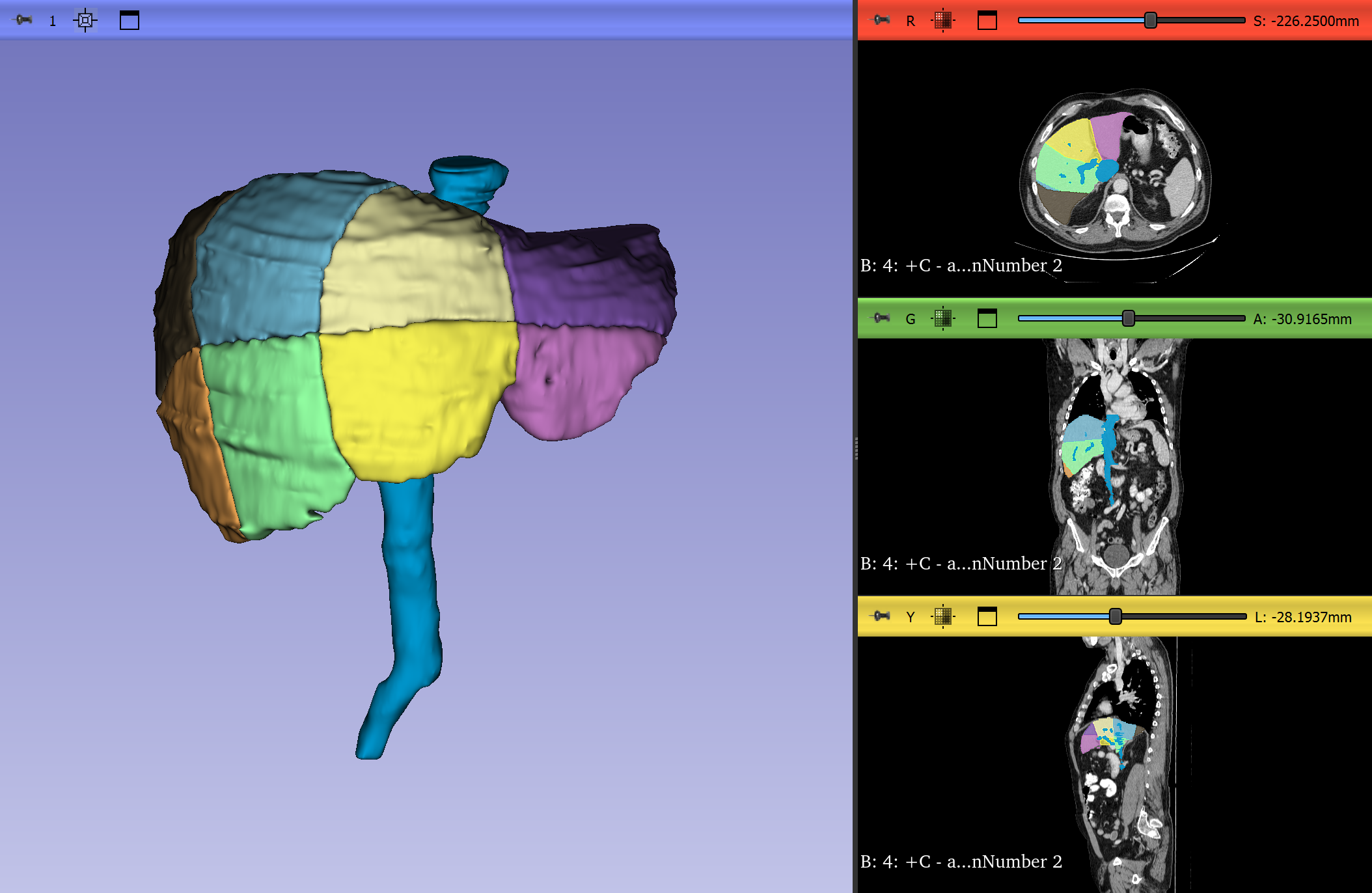}
    \caption{Three-dimensional segmentation of the vascular structure and the eight liver segments in 3D Slicer. The different regions are retained as independent masks to allow the subsequent assignment of spatially differentiated transport properties.}
    \label{fig:liver_segmentation}
\end{figure}

\begin{figure}[h!]
    \centering    \includegraphics[width=0.99\textwidth]{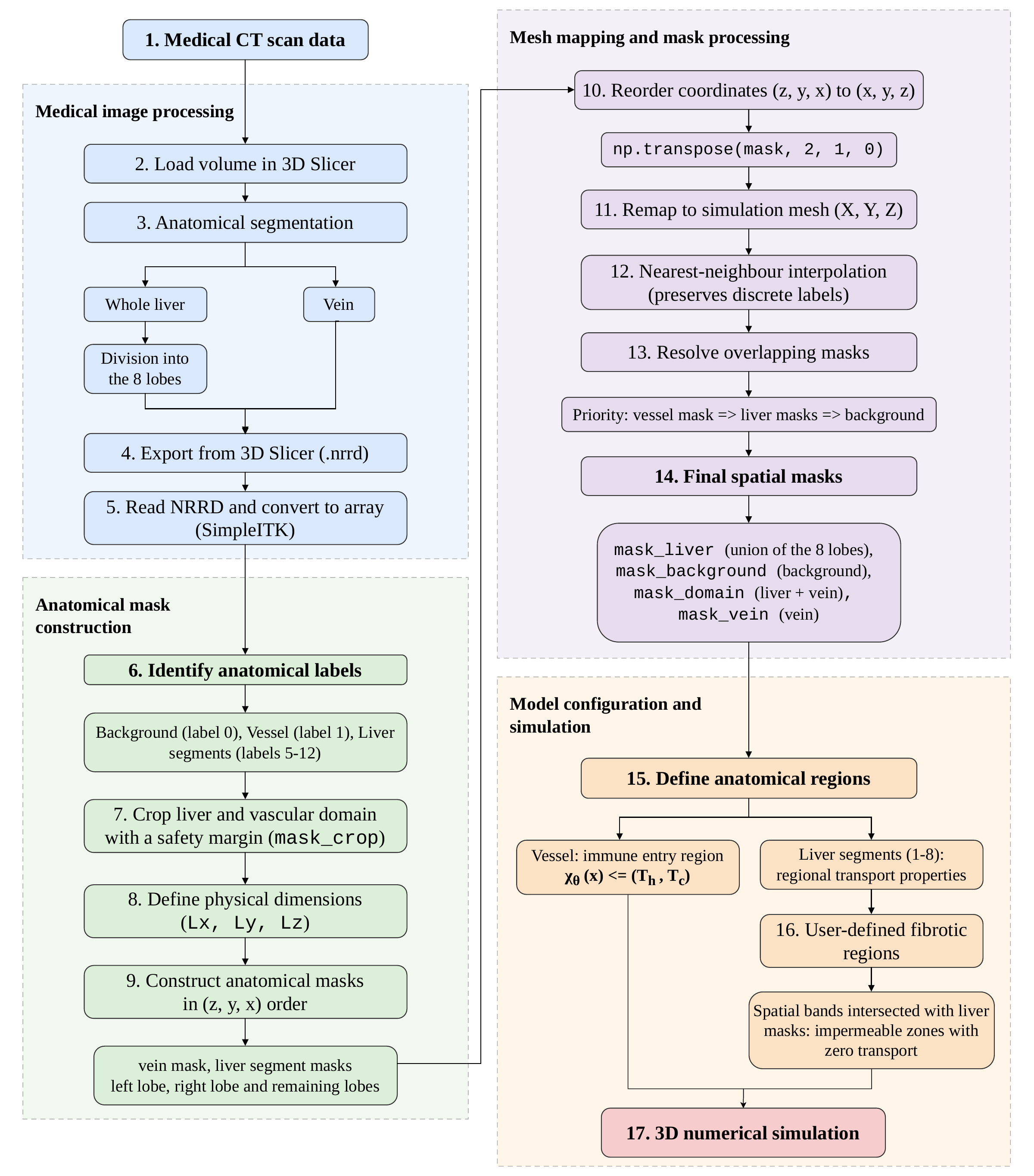}
    \caption{Workflow for constructing the heterogeneous three-dimensional domain from computed tomography images..}
    \label{fig:diagrama_flujo}
\end{figure}

The liver was manually segmented and divided into eight anatomical regions, while the vascular structure was segmented separately to define the localized entry region $\Theta $ 
for $T_h$ and $T_c$ lymphocytes. The resulting anatomical regions were preserved as independent masks, as shown in Figure \ref{fig:liver_segmentation}. The complete workflow, from the initial CT data to the final computational domain, is summarized in Figure \ref{fig:diagrama_flujo}.

The resulting segmentation was exported in NRRD format and subsequently processed with SimpleITK \cite{Lowekamp2013}. The original volume of $(284,512,512)$ voxels was cropped to retain only the region containing the liver and vascular structure, resulting in a volume of $(126,248,293)$ voxels. The physical dimensions of the cropped domain, obtained from the voxel spacing, were
\[
L_x=0.2644\ \mathrm{m}, \qquad
L_y=0.2238\ \mathrm{m}, \qquad
L_z=0.3150\ \mathrm{m}.
\]

The segmentation was then adapted to the numerical mesh of $64\times64\times36$ nodes using nearest-neighbour interpolation, preserving the discrete anatomical labels. Since the arrays obtained from SimpleITK are ordered as $(z,y,x)$, their orientation was converted to the $(x,y,z)$ convention used by the numerical solver. Overlapping regions were subsequently resolved by prioritising the vascular mask over the liver masks.
The resulting masks define the liver, vascular structure and background regions of the computational domain. Spatially heterogeneous diffusion and advection coefficients can therefore be assigned according to the anatomical region. Fibrotic regions are introduced subsequently as spatial bands intersected with the corresponding liver masks; within these regions, the model variables and transport coefficients are set to zero, thereby generating impermeable zones. The resulting heterogeneous domain is then used in the three-dimensional simulations described in Section \ref{sec:3d}.

\subsection{Three-dimensional simulations of a segmented liver}
\label{sec:3d}

The model was applied to a segmented three-dimensional liver geometry comprising the vascular structure and eight liver segments, see Section \ref{subsec:liver_geometry}. 
To analyse the influence of internal transport restrictions on infection evolution, two configurations were compared: one without internal obstacles and another containing fibrotic regions in specific liver segments.

The fibrotic regions were defined using spatial masks restricted to the corresponding segments, with both the model variables and the diffusion and advection coefficients set to zero inside them, see Section \ref{sec:barrier}. Their locations are given in Table \ref{tab:definicion_fibrosis_3D}.
\begin{table}[h!]
\centering
\caption{Geometric definition of the fibrotic regions introduced into the segmented three-dimensional model. Each barrier is initially defined as a spatial band and subsequently restricted to the indicated liver regions.}
\label{tab:definicion_fibrosis_3D}
\resizebox{\textwidth}{!}{
\begin{tabular}{l c c c c}
\hline
Affected lobes & $x$ interval (m) & $y$ interval (m) & $z$ interval (m) \\
\hline
L5--L6 & $0.000 < x < 0.2644$ & $0.000 < y < 0.140$ & $0.162 < z < 0.171$ \\
L6--L7 & $0.000 < x < 0.1130$ & $0.175 < y < 0.185$ & $0.000 < z < 0.3150$ \\
L2--L3 & $0.168 < x < 0.176$ & $0.015 < y < 0.2238$ & $0.000 < z < 0.3150$ \\
L2--L3 & $0.192 < x < 0.200$ & $0.000 < y < 0.100$ & $0.000 < z < 0.3150$ \\
\hline
\end{tabular}
}
\end{table}

The analysis was performed over the time interval $t \in [0, 100]$ and on a mesh of $64\times64\times36$ nodes. The actual spatial domain was $(x,y,z) \in [0, L_x] \times [0, L_y] \times [0, L_z]$, with $L_x=0.2644$ m, $L_y=0.2238$ m, and $L_z=0.3150$ m, starting from the initial conditions
\begin{equation}
    q_1(\mathbf{x},0)=0.1, \quad T_h(\mathbf{x},0)=0, \quad T_c(\mathbf{x},0)=0, \quad q_3(\mathbf{x},0)=0.
    \label{eq:condiciones_iniciales_3D}
\end{equation}
using a stability factor $\sigma=0.15$. The inlet region $\Theta$ was defined from the $mask_{in} = mask_{vein}$ mask of the segmented hepatic vein. 
Both $T_h$ and $T_c$ lymphocytes are introduced in this region.

Because of the computational cost associated with the three-dimensional configuration, both simulations were run on a GPU, in accordance with the performance analysis presented in Section \ref{sec:numerical_validation}.

The reaction parameters and transport coefficients assigned to each liver segment were identical in both configurations and are reported in Tables \ref{tab:parametros_reaccion_3D} and \ref{tab:parametros_difusion_adveccion_3D}. The incorporation of the fibrotic regions is therefore the only difference between the two simulations.
\begin{table}[h!]
\centering
\caption{Reaction parameter values used for the three-dimensional case on the segmented liver domain.}
\label{tab:parametros_reaccion_3D}
\resizebox{\textwidth}{!}{
\begin{tabular}{l c c c c c c c c c c c c c}
\hline
 & $a_1$ & $C_1$ & $\epsilon$ & $\kappa$ & $a_5$ & $a_{2h}$ & $C_{Th}$ & $a_{6h}$ & $a_{6c}$ & $C_{Tc}$ & $a_3$ & $a_{nd}$ & $a_{2c}$ \\
\hline
Value   & 1.0 & 1.0 & 0.05 & 0.05 & 0.036 & 14.0 & 8.0 & 0.2 & 0.03 & 15.0 & 0.6 & 1.0 & 8.0 \\
\hline
\end{tabular}}
\end{table}
\FloatBarrier

\begin{table}[h!]
\centering
\caption{Diffusion and advection coefficient values used for the three-dimensional case on the segmented liver domain, by lobe.}
\label{tab:parametros_difusion_adveccion_3D}
\resizebox{\textwidth}{!}{
\begin{tabular}{l c c c c c c c c c c}
\hline
Coefficient
& Background and barriers
& L2
& L8
& L4a
& L3
& L4b
& L7
& L5
& L6 \\
\hline
$d_1^{\mathrm{ctc}}$        & 0.00 & 0.006 & 0.008 & 0.007 & 0.006 & 0.007  & 0.008 & 0.008 & 0.008 \\
$d_{T_h}^{\mathrm{diff}}$ & 0.00 & 0.08 & 0.09 & 0.085 & 0.08 & 0.085 & 0.09 & 0.09 & 0.09 \\
$d_{T_c}^{\mathrm{chem}}$ & 0.00 & 0.05 & 0.07 & 0.06 & 0.05 & 0.06 & 0.07 & 0.07 & 0.07 \\
$d_3^{\mathrm{diff}}$ & 0.00 & 0.04 & 0.05 & 0.045 & 0.04 & 0.045 & 0.05 & 0.05 & 0.05 \\
\hline
\end{tabular}
}
\end{table}
\FloatBarrier

\subsubsection{Viral clearance regime}

In the absence of fibrotic regions, the variables develop spatially heterogeneous distributions during the transient, determined by the liver geometry, the transport coefficients assigned to each segment, and the vascular entry of the immune response. The $T_h$ and $T_c$ lymphocytes are introduced through the vascular structure and progressively redistribute throughout the rest of the domain. In particular, $T_c$ transport is governed by the cytokine $q_3$ concentration gradient, giving rise to a cytotoxic response that spreads through the tissue and progressively acts on the viral load.

At the final simulation time, the viral load $q_1$ has decreased to virtually zero throughout the liver domain, as have the cytokines $q_3$ and the $T_h$ lymphocytes. The $T_c$ lymphocytes retain only a residual distribution, consistent with their more gradual decline after the infection has been eliminated. The system thus reaches a viral-clearance regime.
The simulation required approximately $24.8$ h of GPU computing time.

\subsubsection{Persistent infection regime}

In the presence of fibrotic regions, the immune response initially propagates from the vascular structure in a manner analogous to the previous configuration. However, the regions in which the variables and transport coefficients have been set to zero introduce internal restrictions that modify the spatial redistribution of the different species. Propagation of the cytotoxic response remains governed by the cytokine $q_3$ gradient, but the barriers limit its access to certain tissue regions and produce a spatial distribution different from that obtained without fibrosis.

\begin{figure*}[t]
    \centering

    \begin{subfigure}[t]{0.49\textwidth}
        \centering
        \includegraphics[width=\linewidth]{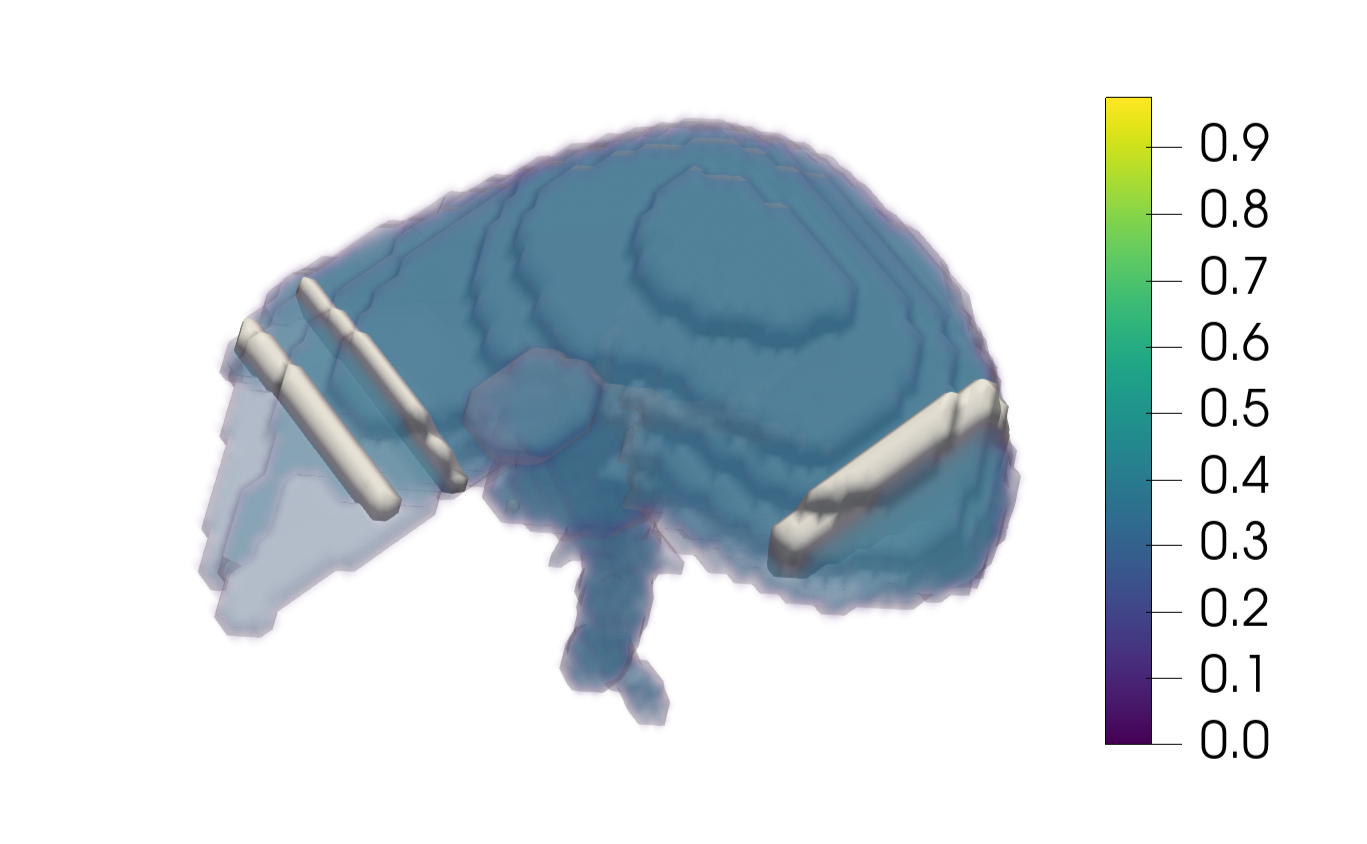}
        \caption{Virus ($q_1$)}
        \label{fig:liver_barrier_q1}
    \end{subfigure}
    \hfill
    \begin{subfigure}[t]{0.49\textwidth}
        \centering
        \includegraphics[width=\linewidth]{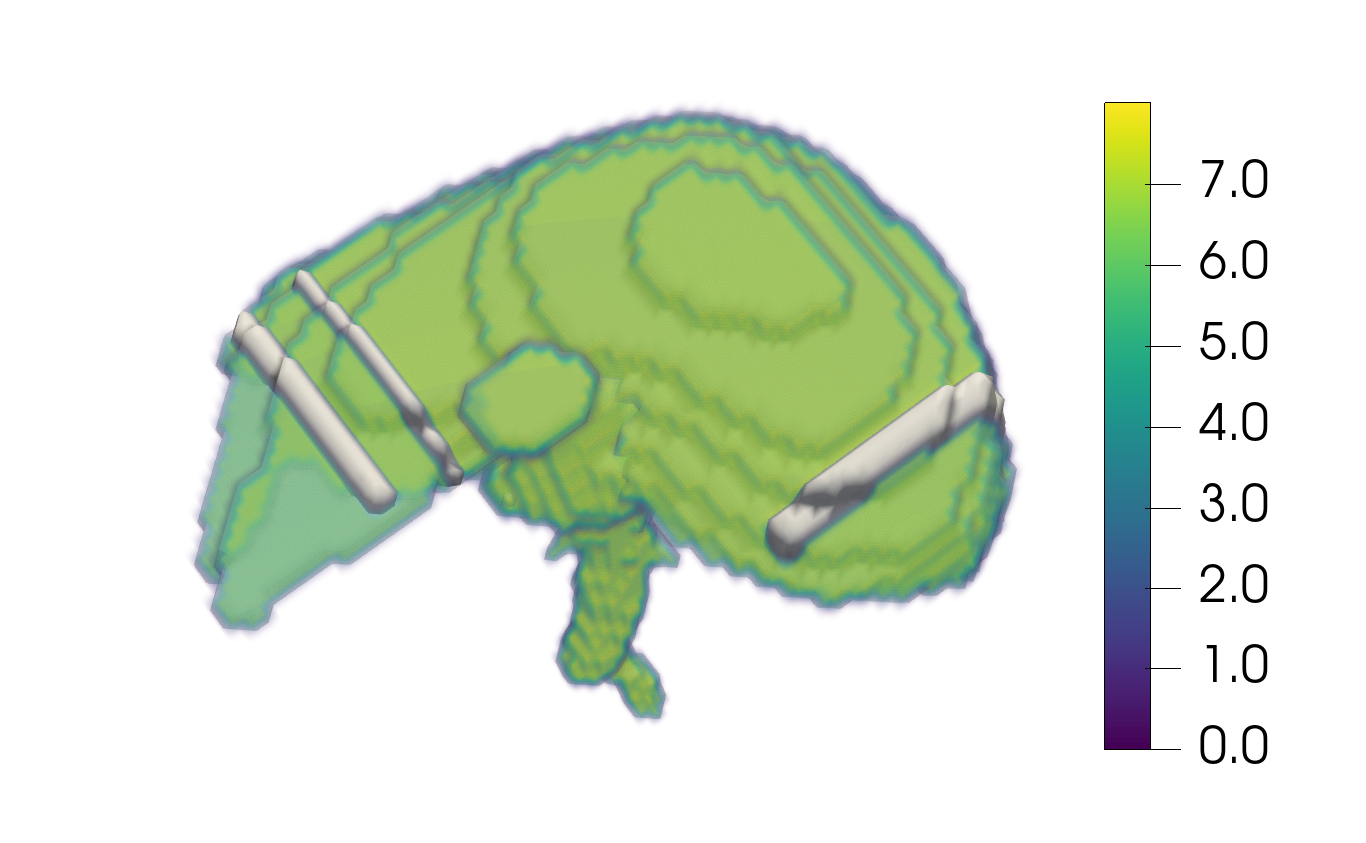}
        \caption{T helper ($T_h$)}
        \label{fig:liver_barrier_Th}
    \end{subfigure}

    \vspace{0.4em}

    \begin{subfigure}[t]{0.49\textwidth}
        \centering
        \includegraphics[width=\linewidth]{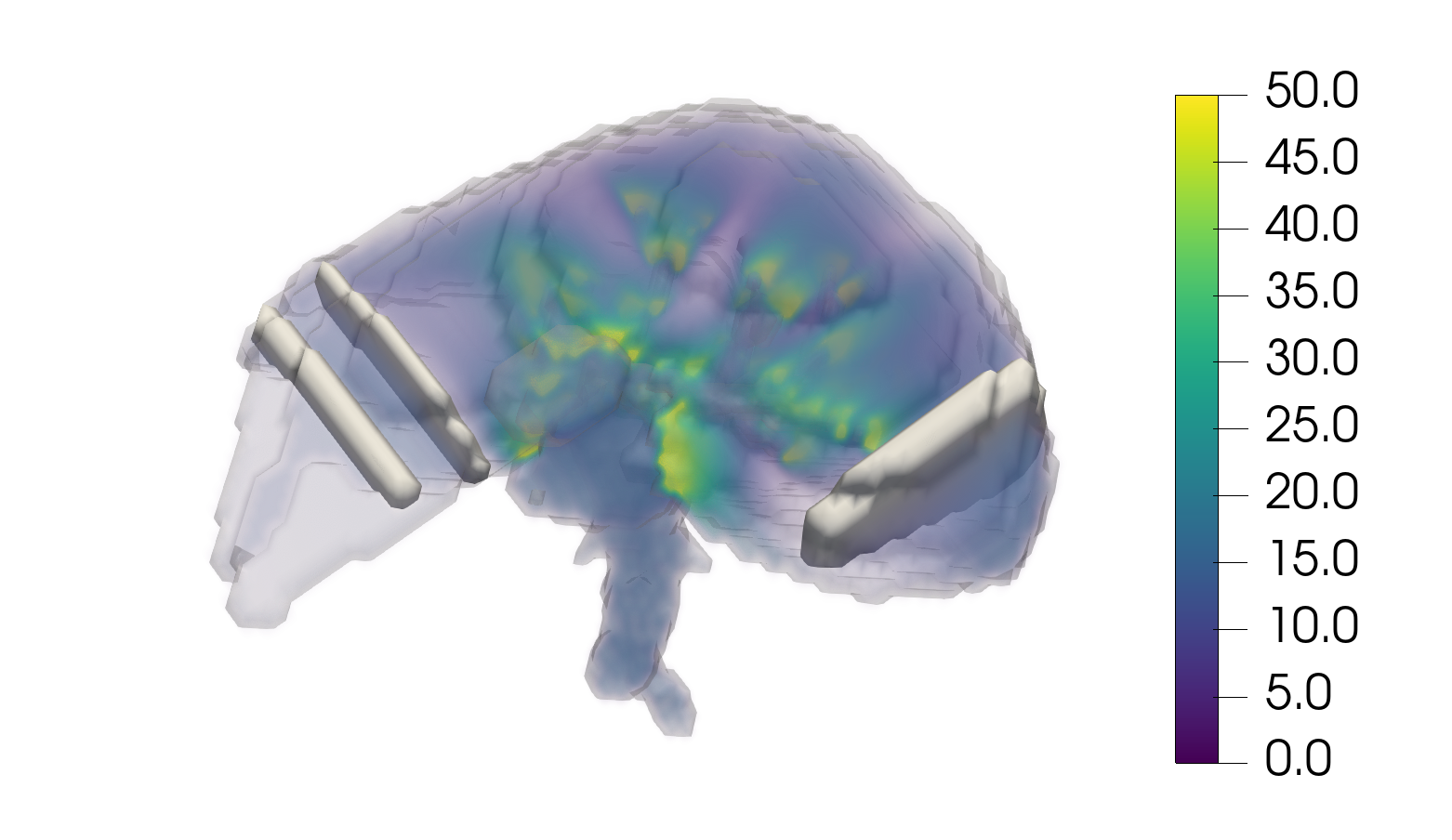}
        \caption{T cytotoxic ($T_c$)}
        \label{fig:liver_barrier_Tc}
    \end{subfigure}
    \hfill
    \begin{subfigure}[t]{0.49\textwidth}
        \centering
        \includegraphics[width=\linewidth]{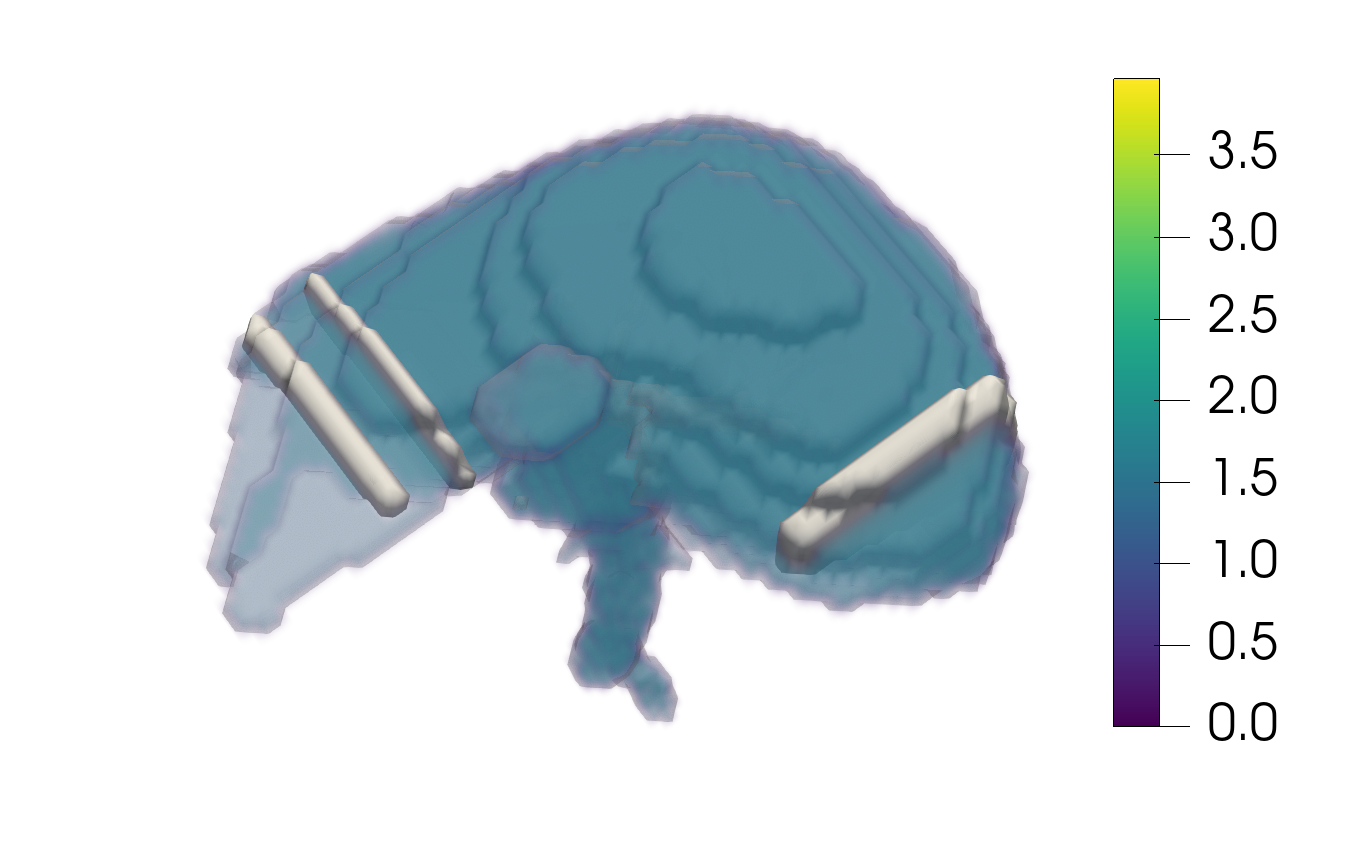}
        \caption{Cytokines ($q_3$)}
        \label{fig:liver_barrier_q3}
    \end{subfigure}

    \caption{Spatial distribution of the four model variables over the
    segmented three-dimensional liver geometry at the end of the simulation,
    $t_f=100$, with fibrotic regions.}
    \label{fig:liver_barrier_variables}
\end{figure*}

This difference is reflected in the final state shown in Figure \ref{fig:liver_barrier_variables}. At $t_f=100$, the viral load $q_1$ remains positive and has a spatially heterogeneous distribution, while the $T_h$ and $T_c$ populations and the cytokines $q_3$ remain active. Unlike the case without fibrosis, the system does not evolve towards viral clearance but towards a persistent-infection regime.
The simulation
required approximately $29.3$ h of GPU computing time.
\FloatBarrier

\begin{figure}[h!]
    \centering
    \begin{minipage}{0.99\textwidth}
        \centering
        \includegraphics[width=\textwidth]{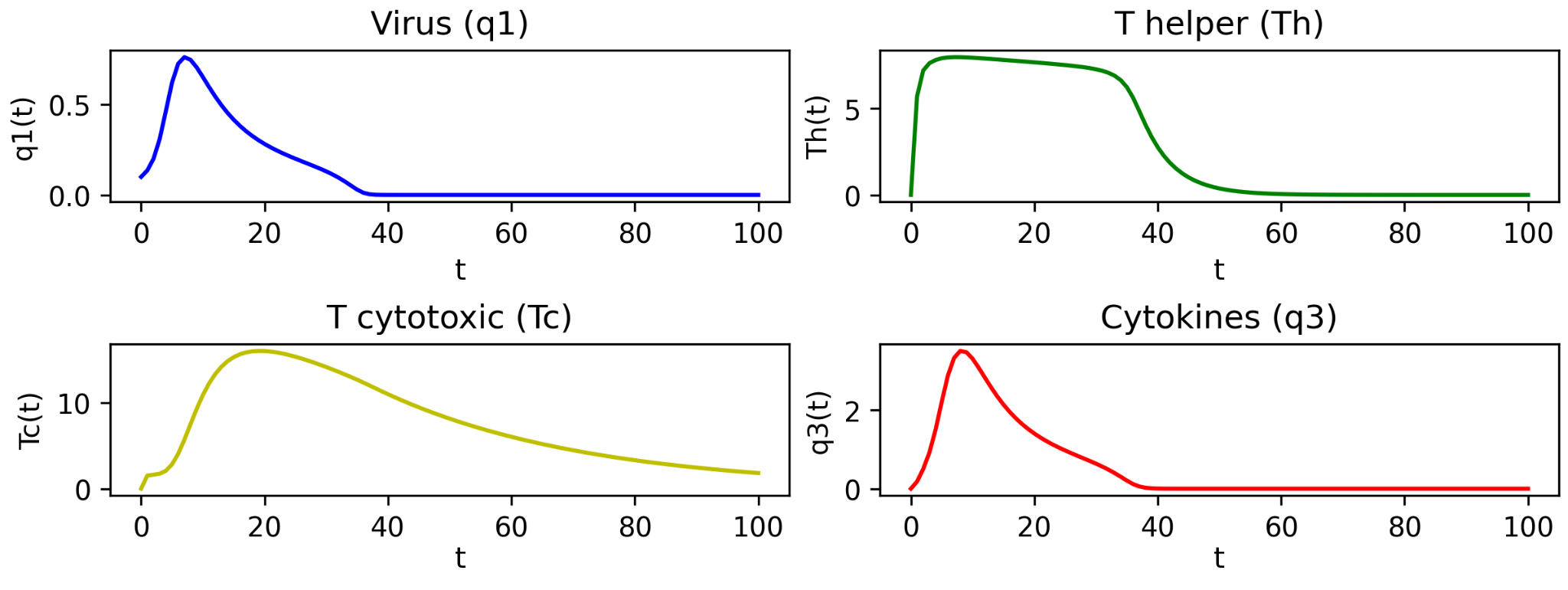}
        \par\smallskip
        (a) Case without fibrosis
    \end{minipage}
    
    \vspace{0.4cm}
    
    \begin{minipage}{0.99\textwidth}
        \centering
        \includegraphics[width=\textwidth]{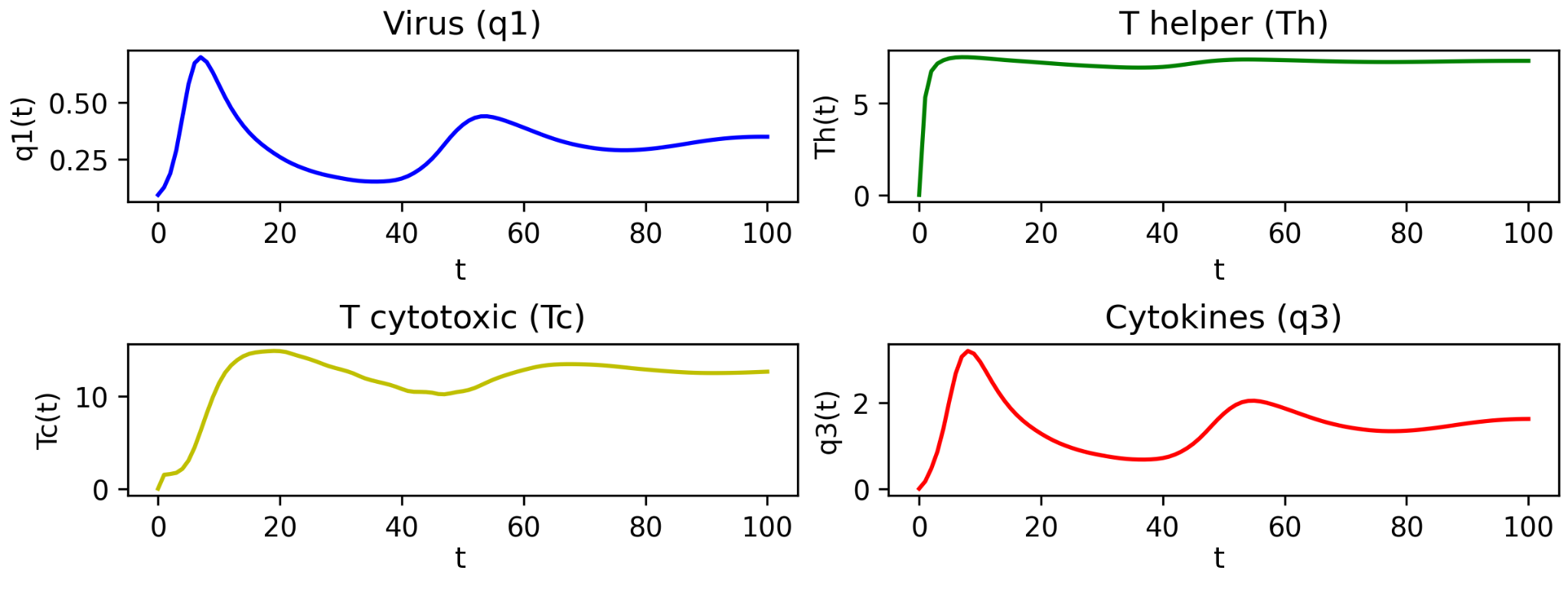}
        \par\smallskip
        (b) Case with fibrosis
    \end{minipage}

    \caption{Temporal evolution of the $L_1$ norm of the four model variables over the segmented three-dimensional liver geometry: (a) configuration without fibrosis and (b) configuration with fibrotic regions.}
    \label{fig:comparacion_3d_temporal}
\end{figure}

The temporal evolution of the two configurations is compared in Figure \ref{fig:comparacion_3d_temporal}. Without fibrosis, after the initial growth of the infection and activation of the immune response, the viral load $q_1$ and cytokines $q_3$ reach a maximum during the early stages of the simulation and subsequently decrease to values close to zero at approximately $t=40$. The $T_h$ lymphocytes show early activation and remain elevated during the first part of the simulation before decreasing rapidly, whereas $T_c$ reaches its maximum later and decays more gradually, retaining residual values at the end of the simulated interval.

In the presence of fibrosis, the evolution exhibits clearly different behaviour. Although the immune response initially reduces the viral load, it does not eliminate it. After reaching an intermediate minimum, $q_1$ increases again and remains positive until the end of the simulation. At the same time, $T_h$, $T_c$, and $q_3$ remain active, confirming the establishment of a persistent-infection regime.
Comparison of the two configurations shows that, with the initial conditions, reaction parameters, and diffusion and advection coefficients assigned to each liver segment held constant, the incorporation of fibrotic regions changes the long-term regime from viral clearance to persistent infection. In the segmented liver geometry, this behaviour reproduces the mechanism previously observed in the two-dimensional domain, where obstruction of spatial transport limits immune-response access to certain regions of the domain and permits local persistence of the viral load.

\section{Conclusions}
\label{sec:conclusions}

In this work, we presented a three-dimensional mathematical model and computational framework for simulating hepatitis B infection dynamics within anatomically realistic liver geometries. The proposed model effectively captures the spatial interaction between virus, T-lymphocyte subpopulations, and cytokines. 
To address the heavy computational burden associated with  3D domains, a parallel finite-difference implementation was developed in Python using NumPy and CuPy, achieving speedup factors exceeding $\times 30$ on GPU architectures compared to traditional CPU execution. 
The resulting framework, released as an open-source library, renders large-scale, high-resolution simulations computationally tractable on segmented CT scan geometries, making personalized computations possible. 
Numerical validation examples demonstrated the model's ability to represent spatial heterogeneities, internal tissue barriers, and localized immune recruitment that govern disease resolution or chronification. 
These spatial heterogeneities are found in reality and may not be described by commonly used models without a spatially resolved three dimensional domain. 
The proposed framework now allows studying the effect of fibrosis during ongoing inflammations. The dynamical development of fibrosis is not yet covered by the model and requires a new system variable that changes over time depending on the local inflammation processes. Overall, this computational tool provides a scalable foundation for future \textit{in-silico} investigations, not only of viral hepatitis, but as well, by replacing the viral dynamics with relevant mechanisms, for Metabolic Dysfunction-Associated Steatohepatitis (MASH), where the inflammation is a consequence of fatty liver tissue.

\appendix

\section{CPU--GPU test case configuration}
\label{app:CPU_GPU_test_case}
The computational performance test presented in Section \ref{subsec: CPU_GPU_performance} was carried out using the same
three-dimensional test case for both the CPU and GPU implementations. The computational domain was defined as $(x,y,z)\in[0,L_x]\times[0,L_y]\times[0,L_z]$, with
$L_x=L_y=L_z=1$. Five uniform meshes were considered, with $n_x=n_y=n_z\in\{50,100,150,200,250\}$.

The system was integrated up to $t_f=2\times10^{-4}$. The time-step
values corresponding to each mesh resolution are reported in Table \ref{tab:comparacion_cpu_gpu}. The output interval was set to $t_d=5$.
All simulations shared the same initial conditions,
\begin{equation}
    q_1(\mathbf{x},0)=1, \quad
    T_h(\mathbf{x},0)=0, \quad
    T_c(\mathbf{x},0)=0, \quad
    q_3(\mathbf{x},0)=0.
    \label{eq:initial_conditions_cpu_gpu}
\end{equation}

The immune-response entry region $\Theta$ was defined as the portion of
a sphere centred at the origin and contained within the computational
domain, with radius
$R_{\mathrm{interior}}=0.21\min(L_x,L_y,L_z)$. In addition, an internal
rectangular barrier was introduced in the region
$0.5\leq x\leq0.6$, $0\leq y\leq0.8$, and $0\leq z\leq1$. The four
model variables were set to zero within this barrier.

The reaction parameters were those reported in Table
\ref{tab:parametros_locales_modelo}, except for $a_5$ and $a_3$, which were set to $a_5=0.06$ and $a_3=0.5$, respectively. The transport
coefficients were those reported in Table
\ref{tab:parametros_distribuidos_modelo}. The same parameter values were used in the CPU and
GPU simulations, so that the computing device was the only difference
between each pair of runs.

For each mesh resolution, the simulation was first executed using the
GPU implementation and subsequently using the CPU implementation. The
reported execution times correspond exclusively to the numerical
computation performed by the solver; the time required to write the
output files was excluded from the measurements.

\section{Visualization of the 3D simulation result}
\label{anexo:Case_3D_github}
To complement the results presented in Section \ref{sec:3d}, animations derived from the three-dimensional simulations are available in the project repository \textit{\href{https://github.com/evaa04/HepGPU-Results}{HepGPU-Results}}. The repository contains results for the various configurations studied throughout the article in \texttt{.gif} format, allowing for the visualization of the temporal and spatial evolution of the model variables across the different case studies.

\section*{Acknowledgements}
ANM and EC acknowledges the Fundación Antonio Gargallo grant 2025/B004. CR acknowledges support through the Novo Nordisk Foundation grant NNF25OC0105537.

\section*{Declaration of interest} none


\begin{thebibliography}{00}

\bibitem{Usuda2024}
D.~Usuda, et~al., Current perspectives of viral hepatitis, World J. Gastroenterol. 30~(18) (2024) 2402--2417. doi: \href{https://doi.org/10.3748/wjg.v30.i18.2402}{10.3748/wjg.v30.i18.2402}.

\bibitem{Reisch2024}
C.~Reisch, S.~Nickel, H.-M.~Tautenhahn, Building up a model family for inflammations, J. Math. Biol. 89~(3) (2024) 29. doi: \href{https://doi.org/10.1007/s00285-024-02126-4}{10.1007/s00285-024-02126-4}.

\bibitem{Kerl2012}
H.-J.~Kerl, D.~Langemann, A.~Vollrath, Reaction diffusion equations and the chronification of liver infections, Math. Comput. Simul. 82~(11) (2012) 2145--2156. doi: \href{https://doi.org/10.1016/j.matcom.2012.04.011}{10.1016/j.matcom.2012.04.011}.

\bibitem{Reisch2019}
C.~Reisch, D.~Langemann, Chemotactic effects in reaction-diffusion equations for inflammation, J. Biol. Phys. 45~(3) (2019) 253--273. doi: \href{https://doi.org/10.1007/s10867-019-09527-3}{10.1007/s10867-019-09527-3}.

\bibitem{Foko2026}
S.~Foko, S.A.~Mohammadi, Complex dynamics in a delayed diffusive cell population model of hepatitis B with capsids, mitotic transmission, and non-cytopathic antiviral mechanisms, Discrete Contin. Dyn. Syst. Ser. B 31 (2026) 94--145. doi: \href{https://doi.org/10.3934/dcdsb.2025104}{10.3934/dcdsb.2025104}.

\bibitem{BoivinChampeaux2026}
C.~Boivin-Champeaux, S.~Schmidt, S.~Balsitis, F.J.~Azeredo, J.S.~Feigelman, A Quantitative Systems Pharmacology (QSP) Model of Acute Hepatitis B Virus Infection: Mechanistic Insights and Foundations for Future Extensions, CPT: Pharmacometrics \& Systems Pharmacology 15 (2) (2026) e70172. doi: \href{https://doi.org}{10.1002/psp4.70172}.


\bibitem{Taylor2008}
Z.A.~Taylor, M.~Cheng, S.~Ourselin, High-speed nonlinear finite element analysis for surgical simulation using graphics processing units, IEEE Trans. Med. Imaging 27~(5) (2008) 650--663. doi: \href{https://doi.org/10.1109/TMI.2007.913112}{\nolinkurl{10.1109/TMI.2007.913112}}.

\bibitem{Madhikar2018}
P.~Madhikar, J.~{\AA}str{\"o}m, J.~Westerholm, M.~Karttunen, CellSim3D: GPU accelerated software for simulations of cellular growth and division in three dimensions, Comput. Phys. Commun. 232 (2018) 206--213. doi: \href{https://doi.org/10.1016/j.cpc.2018.05.024}{10.1016/j.cpc.2018.05.024}.

\bibitem{Hesam2021}
A.~Hesam, L.~Breitwieser, F.~Rademakers, Z.~Al-Ars, GPU acceleration of 3D agent-based biological simulations, in: 2021 IEEE International Parallel and Distributed Processing Symposium Workshops (IPDPSW), 2021, pp.~210--217. doi: \href{https://doi.org/10.1109/IPDPSW52791.2021.00040}{10.1109/IPDPSW52791.2021.00040}.

\bibitem{Hirsch2007}
C.~Hirsch, Numerical Computation of Internal and External Flows: The Fundamentals of Computational Fluid Dynamics, 2nd ed., Butterworth--Heinemann, Oxford, 2007. doi: \href{https://doi.org/10.1016/B978-0-7506-6594-0.X5037-1}{10.1016/B978-0-7506-6594-0.X5037-1}.

\bibitem{Fedorov2012}
A.~Fedorov, R.~Beichel, J.~Kalpathy-Cramer, et~al., 3D Slicer as an image computing platform for the Quantitative Imaging Network, Magn. Reson. Imaging 30~(9) (2012) 1323--1341. doi: \href{https://doi.org/10.1016/j.mri.2012.05.001}{10.1016/j.mri.2012.05.001}.

\bibitem{Rister2020}
B.~Rister, D.~Yi, K.~Shivakumar, et~al., CT-ORG, a new dataset for multiple organ segmentation in computed tomography, Sci. Data 7~(1) (2020) 381. doi: \href{https://doi.org}{10.1038/s41597-020-00715-8}.

\bibitem{Kanel2017}
G.C.~Kanel, Pathology of Liver Diseases, John Wiley \& Sons, Hoboken, NJ, 2017. doi: \href{https://doi.org/10.1002/9781118894903}{10.1002/9781118894903}.

\bibitem{Lowekamp2013}
B.C.~Lowekamp, D.T.~Chen, L.~Ib{\'a}{\~n}ez, D.~Blezek, The design of SimpleITK, Front. Neuroinform. 7 (2013) 45. doi: \href{https://doi.org/10.3389/fninf.2013.00045}{10.3389/fninf.2013.00045}.

\bibitem{Okuta2017}
R.~Okuta, Y.~Unno, D.~Nishino, S.~Hido, C.~Loomis, CuPy: A NumPy-compatible library for NVIDIA GPU calculations, in: Proceedings of the Workshop on Machine Learning Systems (LearningSys), Thirty-first Annual Conference on Neural Information Processing Systems, 2017. \url{https://learningsys.org/nips17/assets/papers/paper_16.pdf} (accessed 20 August 2026).

\bibitem{Harris2020}
C.R.~Harris, K.J.~Millman, S.J.~van der Walt, et~al., Array programming with NumPy, Nature 585~(7825) (2020) 357--362. doi: \href{https://doi.org/10.1038/s41586-020-2649-2}{10.1038/s41586-020-2649-2}.

\bibitem{Zarr2025}
Zarr Developers, zarr-developers/zarr-python: v2.18.5, Zenodo, 2025. doi: \href{https://doi.org/10.5281/zenodo.15101732}{10.5281/zenodo.15101732}.

\bibitem{Ahrens2005}
J.~Ahrens, B.~Geveci, C.~Law, ParaView: An end-user tool for large data visualization, in: C.D.~Hansen, C.R.~Johnson (Eds.), The Visualization Handbook, Elsevier, Burlington, MA, 2005, pp.~717--731. doi: \href{https://doi.org/10.1016/B978-012387582-2/50038-1}{10.1016/B978-012387582-2/50038-1}.

\bibitem{Wasserthal2023}
J.~Wasserthal, H.-C.~Breit, M.T.~Meyer, et~al., TotalSegmentator: Robust segmentation of 104 anatomic structures in CT images, Radiol. Artif. Intell. 5~(5) (2023) e230024. doi: \href{https://doi.org/10.1148/ryai.230024}{10.1148/ryai.230024}.

\bibitem{Navas2020}
A.~Navas Montilla, P.~Garc{\'i}a Navarro, J.~Fern{\'a}ndez Pato, Introducci{\'o}n a la mec{\'a}nica de fluidos computacional, Prensas de la Universidad de Zaragoza, Zaragoza, 2020, ISBN 978-84-1340-233-8.

\bibitem{Jiang2015}
B.~Jiang, W.~Dai, A.~Khaliq, M.~Carey, X.~Zhou, L.~Zhang, Novel 3D GPU based numerical parallel diffusion algorithms in cylindrical coordinates for health care simulation, Math. Comput. Simul. 109 (2015) 1--19. doi: \href{https://doi.org/10.1016/j.matcom.2014.07.003}{10.1016/j.matcom.2014.07.003}.

\bibitem{Hunter2007}
J.D.~Hunter, Matplotlib: A 2D graphics environment, Comput. Sci. Eng. 9~(3) (2007) 90--95. doi: \href{https://doi.org/10.1109/MCSE.2007.55}{10.1109/MCSE.2007.55}.

\bibitem{Schroeder2006}
W.~Schroeder, K.~Martin, B.~Lorensen, The Visualization Toolkit: An Object-Oriented Approach to 3D Graphics, 4th ed., Kitware, Clifton Park, NY, 2006, ISBN 978-1-930934-19-1.


\end{thebibliography}
\end{document}